\documentclass[10pt]{article}
\usepackage[all]{xy}
\usepackage{amsmath}
\usepackage{amsfonts}
\usepackage{amssymb}
\usepackage{amscd}
\usepackage{amsthm}
\usepackage{latexsym}
\usepackage{amsbsy}
\usepackage{color}
\newtheorem{lemma}{Lemma}[section]
\newtheorem{proposition}{Proposition}[section]
\newtheorem{theorem}{Theorem}[section]
\newtheorem{corollary}{Corollary}[section]
\newtheorem{definition}{Definition}[section]
\newtheorem{example}{Example}[section]

\newtheorem{remark}{Remark}[section]

\newcommand{\nc}[2]{\newcommand{#1}{#2}}
\newcommand{\rnc}[2]{\renewcommand{#1}{#2}}
\rnc{\theequation}{\thesection.\arabic{equation}}
\nc{\beq}{\begin{equation}}
\nc{\eeq}{\end{equation}}
\rnc{\[}{\beq}
\rnc{\]}{\eeq}
\nc{\ba}{\begin{array}}
\nc{\ea}{\end{array}}
\nc{\bea}{\begin{eqnarray}}
\nc{\beas}{\begin{eqnarray*}}
\nc{\eeas}{\end{eqnarray*}}
\nc{\eea}{\end{eqnarray}}
\nc{\be}{\begin{enumerate}}
\nc{\ee}{\end{enumerate}}
\nc{\bd}{\begin{diagram}}
\nc{\ed}{\end{diagram}}
\nc{\bi}{\begin{itemize}}
\nc{\ei}{\end{itemize}}
\nc{\bpr}{\begin{proposition}}
\nc{\bth}{\begin{theorem}}
\nc{\ble}{\begin{lemma}}
\nc{\bco}{\begin{corollary}}
\nc{\bre}{\begin{remark}}
\nc{\bex}{\begin{example}}
\nc{\bde}{\begin{definition}}
\nc{\ede}{\end{definition}}
\nc{\epr}{\end{proposition}}
\nc{\ethe}{\end{theorem}}
\nc{\ele}{\end{lemma}}
\nc{\eco}{\end{corollary}}
\nc{\ere}{\end{remark}}
\nc{\eex}{\end{example}}
\nc{\bpf}{{ Proof. ~~}}
\nc{\epf}{\hfill\mbox{$\square$}\vspace*{3mm}}
\nc{\hsp}{\hspace*}
\nc{\vsp}{\vspace*}

\nc{\mc}{\mathcal{C}}
\nc{\ws}{\widetilde{S}}
\nc{\wdel}{\widetilde{\delta}}
\nc{\wsigma}{\widetilde{\sigma}}
\nc{\wtau}{\widetilde{\tau}}
\nc{\wvfy}{\widetilde{\varphi}}

\nc{\ra}{\rightarrow} \nc{\lra}{\longrightarrow}
\nc{\lla}{\longleftarrow}

\nc{\epn}{\varepsilon}
\nc{\si}{\psi}
\nc{\Del}{\Delta}
\nc{\del}{\delta}
\nc{\ro}{\rho}
\nc{\fy}{\phi}
\nc{\tu}{\tau}
\nc{\vfy}{\varphi}
\nc{\sig}{\sigma}

\nc{\alf}{\alpha} \nc{\bt}{\beta} \nc{\mh}{\mathcal{H}}
\nc{\mbc}{\mathbb{C}} \nc{\mfg}{\mathfrak{g}} \nc{\ug}{U(\mfg)}
\nc{\tg}{T(\mfg)}

\nc{\sg}{S(\mfg)}

 \nc{\opl}{\oplus} \nc{\bopl}{\bigoplus}
\nc{\Lam}{\Lambda} \nc{\lam}{\lambda} \nc{\weg}{\wedge}
\nc{\bweg}{\bigwedge} \nc{\blt}{\bullet}
\nc{\ydh}{{}_{H}^{H}\mathcal{YD}}

\nc{\nn}{\nonumber}

\nc{\la}{\leftarrow}
\nc{\da}{\downarrow}
\nc{\ua}{\uparrow}
\nc{\ot}{\otimes}

\nc{\uh}{\underline{H}}
 \nc{\udel}{\underline{\Del}}
\nc{\us}{\underline{S}}
 \nc{\um}{\underline{m}}
\nc{\uepn}{\underline{\epn}}
\nc{\ueta}{\underline{\eta}}

\begin{document}
\title{Cup Coproducts in Hopf Cyclic Cohomology }
\author {Mohammad Hassanzadeh and  Masoud Khalkhali \\ Department of Mathematics,  University of Western Ontario \\ London, Ontario, Canada}
\date{}
\maketitle
\begin{abstract}
We define cup coproducts for Hopf cyclic cohomology  of Hopf algebras and for its dual theory. We show that for universal enveloping algebras and group algebras our coproduct recovers
the standard coproducts on Lie algebra homology and group homology, respectively.
\end{abstract}

\section{ Introduction}
In this paper we define cup coproducts for Hopf cyclic cohomology of Hopf algebras by establishing K\"{u}nneth formulas for Hopf cyclic cohomology. Cup product  for  cyclic cohomology  of algebras was first defined by Connes in   \cite{c2}. K\"{u}nneth formulas and (co)products  for  cyclic theory of algebras has been established  by several authors \cite{bav, bn, em,  ka, ka1, kas1} (cf. also Loday's book \cite{L} and references therein for a full account).  Ever since the discovery of Hopf cyclic cohomology by Connes and Moscovici and with computations carried on in \cite{cm4, cm3},
 it was clear that this theory and its dual counterpart \cite{kr1}   is  a noncommutative analogue of Lie algebra homology and group homology. This of course poses a natural question:
  to what extent cup products in Lie algebra and group cohomology can be extended to Hopf cyclic cohomology. This question is answered in this paper.  We use the theory of cyclic modules and the definition of Hopf cyclic cohomology (with coefficients) in terms of cyclic modules \cite{cm4,hkrs2} to define our coproducts in a natural way via Eilenberg-Zilber isomorphisms for cyclic modules. This gives us an external coproduct. By specifying to cocommutative Hopf algebras we then define cup  coproducts for Hopf cyclic (co)homology
  of cocommutative Hopf algebras.

This paper is organized as follows. In Sections 2 and  3  we recall the Eilenberg-Zilber isomorphism and K\"{u}nneth formula   for cyclic homology. In Section 4 we recall the definition of
 Hopf cyclic cohomology with coefficients in an stable anti-Yetter Drinfeld module and work out  the Eilenberg-Zilber isomorphism and K\"{u}nneth formula in the context of Hopf cyclic cohomology. In Section 5  we define external coproducts for Hopf cyclic cohomology and in Section 6 by restricting to cocommutative Hopf algebras we obtain internal coproducts on Hopf cyclic
 cohomology.     It is known that when restricted to universal enveloping algebras,  Hopf cyclic cohomology reduces to Lie algebra homology \cite{cm4}. In Sections 6 we show  that under these isomorphisms our coproduct  coincides with the coproduct   in Lie algebra homology. In Section 7 we look at the dual Hopf cyclic homology as defined in \cite{kr1} and define a coproduct on it. It is shown in
\cite{kr1} that when restricted to group  algebras, Hopf cyclic homology is isomorphic  to group homology. In this section we check that under this isomorphism our coproduct coincides with coproduct in  group homology.

It should be mentioned that the term cup product is also used for  a  different type of operation  in Hopf cyclic theory. Thus the cup products defined and studied in
\cite{gor, hkrs2,  kay1, kay, kr3, Ran} involve an action
  of a Hopf algebra on  a (co)algebra and  is a pairing between  Hopf cyclic cohomology of the given (co)algebra and Hopf cyclic cohomology of the acting Hopf algebra.  This operation is in fact
  an extension of Connes-Moscovicvi's
  characteristic map from \cite{cm3}.


\section{ Eilenberg-Zilber isomorphisms }
 In this section we recall the Eilenberg-Zilber isomorphisms for Hochschild, cyclic and periodic cyclic (co)homology of (co)cyclic modules where they will be used crucially in future sections \cite{bav, gj, khra}. Given a (co)cyclic module $(C, \delta, s , \tau)$ \cite{cm2, L}, where $\delta$, $s$ and $\tau$ are (co)faces, (co)degeneracies and (co)cyclic maps respectively, we denote its Hochschild differential in cohomology by $b^n=\sum_{i=0}^{n} (-1)^{i}\delta_{i}$, and Connes' differential by $B^n=(\sum_{i=0}^{n} (-1)^{ni}\tau_{n}^i)s_n\tau_{n}(1-\tau_n)$ . Any (co)cyclic module defines a mixed complex $(C, b, B)$ in a natural way \cite{cm2,kas,L}. Suppose $(C^{n},\delta^{n}_{i},s^{n}_{i},\tau^{n})$ and $(C'^{n},\delta'^{n}_{i},s'^{n}_{i},\tau'^{n})$ are two cocyclic modules.  The diagonal $C\times C'$  is the cocyclic module  $ ((C\times C')^{n},\delta^{n}_{i}\otimes \delta'^{n}_i,s^{n}_i\otimes s'^{n}_i,\tau_{n}\otimes \tau'_{n}) $
where $$  (C\times C')^{n} = C^{n}\otimes C'^{n}. $$ The tensor product complex  $ C\otimes C'$ given by $$(C\otimes C')^{n} = \bigoplus_{{p+q=n}} C^{p}\otimes C'^{q},$$
  is not a cocyclic module but it has a  mixed complex structure given by  $$ b_{n} = \bigoplus_{{i+j=n}} ( b_{i}\otimes id_{j} + (-1)^{j}id_{i}\otimes b'_{j}),$$
 and similarly for $B$.

 We have two natural maps between these complexes. The Alexander-Whitney map $$AW_{n}:(C\otimes C')^{n}  \longrightarrow (C\times C')^{n} ,$$
 is given by
 \[ AW_{n} = \bigoplus_{{i+j=n}} AW_{i,j},\]where $$AW_{i,j} : C^{i}\otimes C'^{j}\longrightarrow C^{n}\otimes C'^{n},  $$ are defined by
  \[AW_{i,j} = \delta^{n}_{n}\delta^{n-1}_{n-1} \cdots \delta^{i+1}_{i+1}\otimes \delta^{n}_{0}\delta^{n-1}_{0}\cdots \delta^{i+1}_{0}. \]
The Shuffle map
 $$sh_{n}: (C\times C')^{n}\longrightarrow (C\otimes C')^{n},$$ is given by
 \[ sh_{n}= \bigoplus_{{i+j=n}}sh_{i,j}, \] where $$sh_{i,j} : (C\times C')^{n} \longrightarrow C^{i}\otimes C'^{ j},$$ are defined by
\[sh_{i,j} = \sum_{{\sigma}} sign(\sigma)s_{\sigma(i+1)}^{i}\cdots s_{\sigma(n)}^{n-1}\otimes s_{\sigma(1)}^{j}\cdots s_{\sigma(i)}^{n-1},\] and $\sigma$ runs over all $(i,j)$-shuffles. By an $(i,j)$ shuffle we mean a permutation $ \sigma $ on elements  $\{1,2,\dots,n\}$  such that
$$ \sigma(1)<\sigma(2)<\cdots <\sigma(i)    , \qquad   \sigma(i+1)<\sigma(i+2)<\cdots<\sigma(n).$$One has  $[b,sh]=0$ and $[b,AW]=0$ in the normalized setting \cite{cm2, L}. It can be shown that
$$sh\circ AW=id    , \qquad   AW\circ sh = bh + hb + id,$$ for some chain homotopy map $h$ \cite{L}. So they induce inverse isomorphisms on Hochschild  cohomology, called the Eilenberg-Zilber isomorphism:
  \[ HH^{n}(C\otimes C')\cong HH^{n}(C\times C')     \qquad  \forall n\geq 0.\label{EZIH} \]

  In the dual case of cyclic modules, the map

  $$ AW_{n} : (C\times C')_{n}\longrightarrow (C\otimes C')_{n}, $$is given by
 \[ AW_{n} = \bigoplus_{{i+j=n}} AW_{i,j},\] where $ AW_{i,j} : C_{n}\otimes C'_{n}\longrightarrow C_{i}\otimes C'_{j}$ are given by

 \[AW_{i,j} = \delta^{i+1}_{i+1}\delta^{i+2}_{i+2} \cdots \delta^{n}_{n}\otimes \delta^{i+1}_{0}\delta^{i+2}_{0}\cdots \delta^{n}_{0}. \]
Also
 $$ sh_{n}:(C\otimes C')_{n}\longrightarrow (C\times C')_{n}, $$ is given by
 \[ sh_{n}= \bigoplus_{{i+j=n}}sh_{i,j},  \] where $sh_{i,j}: C_{i}\otimes C'_{j}\longrightarrow C_{n}\otimes C'_{n}$ can be defined as:

 \[sh_{i,j} = \sum_{{\sigma}} sign(\sigma)s_{\sigma(n)}^{n-1}\cdots s_{\sigma(i+1)}^{i}\otimes s_{\sigma(i)}^{n-1}\cdots s_{\sigma(1)}^{j}.\]

 Now we state the Eilenberg-Zilber isomorphism for cyclic cohomology. To this end, we define another map namely the cyclic shuffle map as follows \cite{ku, L}. First we define an $(i,j)$-cyclic shuffle.  Let $ i, j, n \in \mathbb{N}$ with $n=i+j$. Consider the permutation $\sigma$ on the $n$ elements $\{1,\dots,n\}$ obtained by first performing a cyclic permutation $p$ times on $\{1,\dots,i\}$ and a cyclic permutation $q$ times on $\{i+1,\dots,i+j\}$ and there after applying an $(i,j)$-shuffle to the combined result. We call $\sigma$ an $(i,j)$-cyclic shuffle if $1$ appears before $i+1$ in the resulting sequence. One can define another map namely $sh'$ which is in fact a cyclic version of the shuffle map:
 $$sh'_{n}:(C\times C')^{n}\longrightarrow (C\otimes C')^{n-2},$$ where
  \[sh'_{n}=\bigoplus_{{i+j=n}}sh'_{i,j},\]
and  $sh'_{i,j}:C^{n}\otimes C'^{n}\longrightarrow C^{i-1}\otimes C'^{j-1}, $ are given by

 \[sh'_{i,j} =
  \sum_{{\sigma}}sign(\sigma)s^{i-1}_{i-p-1}\tau^{p+1}_{i}s^{i}_{\sigma(i+1)}\cdots s^{n-1}_{\sigma(n)}\otimes s_{j-q-1}^{j-1}\tau_{j}^{q+1}s_{\sigma(1)}^{j}\cdots s_{\sigma(i)}^{n-1}.\]Here $\sigma$ runs over all $(i,j)$-cyclic shuffles.

We denote the $(b,B)$-bicomplex of a (co)cyclic module $C$ by $\mathcal{B}(C)$ and its total complex by $Tot\mathcal{B}(C)$. It is a mixed complex. Recall that an $S$-map \cite{kas, L} between mixed complexes is a map of complexes
$f: Tot\mathcal{B}(C)\longrightarrow Tot\mathcal{B}(C')$ which  commutes with Connes' periodicity map $S$. In fact an $S$-map has a matrix representation:

  \begin{displaymath}
  \left( \begin{array}{cccc}
  f_{0} & f_{1} & f_{2} & \ldots \\
  f_{-1} & f_{0} & f_{1} & \ldots \\
  f_{-2} & f_{-1} & f_{0} & \ldots \\
  \vdots & \vdots & \ddots
  \end{array} \right)
  \end{displaymath}where  $f_{i} : C^{*}\longrightarrow C'^{*-2i},$ and
  \[ [B, f_{i+1}]+ [b, f_i]=0 .\label{eer}\] An $S$-map induces a map on the  level of cyclic cohomology. Now we define the desired $S$-map $\widetilde{Sh}:Tot\mathcal{B}(C\times C')\longrightarrow Tot \mathcal{B}(C\otimes C')$ as:

 \begin{displaymath}
 {\widetilde{Sh}}=
 \left( \begin{array}{cccc}
  sh & sh' & 0 & \ldots \\
  0 & sh & sh' & \ldots \\
  0 & 0 & sh & \ldots \\
  \vdots & \vdots & \ddots
  \end{array} \right)
\end{displaymath} The condition  \eqref{eer} reduces to the following relations in the normalized setting:

(1) $[b,sh]=0,$

(2) $[B,sh]+[b,sh']=0,$

(3) $[B,sh']=0. $

Also there is another $S$-map $\widetilde{AW}:Tot\mathcal{B}(C\otimes C')\longrightarrow Tot\mathcal{B}(C\times C'),$ given by

\begin{displaymath}
 {\widetilde{AW}} =
  \left( \begin{array}{cccc}
  AW & AW' & 0 & \ldots \\
  0 & AW & AW' & \ldots \\
  0 & 0 & AW& \ldots \\
  \vdots & \vdots & \ddots
  \end{array} \right)
\end{displaymath} Here $AW'_n=  \varphi_n B_n AW_n: (C\otimes C')^{n}\longrightarrow (C\times C')^{n-2},$ where \cite{bav, R}

\begin{eqnarray}
& \varphi_n= \sum(-1)^{n-p-q+\sigma(\alpha,\beta)}&(s_{\beta_{q}+n-p-q}\cdots s_{\beta_{1}+n-p-q}s_{n-p-q-1}
 \delta_{n-q+1}\cdots \delta_{n}) \nonumber \\
&&\otimes (s_{\alpha_{p+1}+n-p-q}\cdots s_{\alpha_{1}+n-p-q}\delta_{n-p-q}\cdots
\delta_{n-q-1}).\nonumber
\end{eqnarray}
 The sum is taken over all $0\leq q\leq n-1$, $0\leq p \leq n-q-1,$ and $(\alpha, \beta)$ runs over all $(p+1,q)$-shuffles where  $\sigma(\alpha,\beta)=\sum(\alpha_{i}-i+1).$ One can show \cite{bav}
 $$\widetilde{Sh}\circ \widetilde{AW}=id,   \qquad  \widetilde{AW}\circ \widetilde{Sh}=(b+B)\circ h + h\circ(b+B) + id,$$ for some homotopy map $h:C\otimes C'\longrightarrow C\times C'$. Thus we obtain the Eilenberg-Zilber isomorphism in cyclic cohomology \cite{bav,khra,L}:
\[ HC^{*}(C\times C')\cong HC^{*}(C\otimes C').\label{EZIC}\] In the case of cyclic modules, one can define
\[sh'_{n}:(C\otimes C')_{n}\longrightarrow (C\times C')_{n+2}  \qquad  \forall n\geq 0, \] and
\[AW'_{n}:(C\times C')_{n}\longrightarrow (C\otimes C')_{n+2}  \qquad  \forall n\geq 0. \]

 This enables us to have the $S$-maps $\widetilde{Sh}: Tot\mathcal{B}(C\otimes C')\longrightarrow Tot\mathcal{B}(C\times C'),$  and $\widetilde{AW}: Tot\mathcal{B}(C\times C')\longrightarrow Tot\mathcal{B}(C\otimes C')$  which induce the Eilenberg-Zilber isomorphism in cyclic homology.

We have the Eilenberg-Zilber isomorphism for periodic cyclic cohomology as follows. One knows that $\underrightarrow{\lim}Tot^{2n+ \ast}\mathcal{B}(C)= \oplus_{{n\geq0}}C^{2n+ \ast},$ where $\ast = 0,1$ and  the direct limit is with respect to Connes' periodicity map $S$. Since direct limit commutes with homology we have
$$\underrightarrow{\lim}HC^{2n+\ast}(C)=  \underrightarrow{\lim}H(Tot^{2n+\ast}\mathcal{B}(C))=H(\underrightarrow{\lim}Tot^{2n+\ast}\mathcal{B}(C))=H(\oplus_{{n\geq0}} C^{2n+\ast})=HP^{\ast}(C),$$
where $\ast=0,1$ . Since
$$HP^{\ast}( C\otimes C')= \underrightarrow{\lim}HC^{2n+\ast}(C\otimes C')\cong  \underrightarrow{\lim}HC^{2n+\ast}(C \times C')=HP^{\ast}(C \times C'),$$
  we obtain the Eilenberg-Zilber isomorphism for periodic cyclic cohomology:
\[HP^{*}(C\otimes C')\cong HP^{*}(C\times C'), \label{abcde}\]
  The explicit maps which induce the isomorphism \eqref{abcde}, are infinite matrix versions of $\widetilde{Sh}$ and $\widetilde{AW}$.

    There is an obstacle to get the Eilenberg-Zilber isomorphism for periodic cyclic homology. The problem comes from the fact that homology does not commute with inverse limit in general. The Mittag-Leffler condition \cite{em2} on the inverse system $(HC(C)[-2m],S)_{m}$ guarantees the commutativity of homology and inverse limit. If we have this condition, with a similar argument we obtain $$HP_{n}(C \otimes C')\cong HP_{n}(C\times C').$$


\section{K\"{u}nneth formulas}

In this section we  recall the K\"{u}nneth formula for Hochschild and cyclic (co)homology of (co)cyclic modules and introduce a formula for the periodic case. If $C$ and $C'$ are cocyclic objects in the category of vector spaces,  we have the following  K\"{u}nneth formula:
 \[HH^{n}(C\otimes C'){\cong} \bigoplus_{{p+q=n}} HH^{p}(C)\otimes HH^{q}(C') \qquad  \forall \, n\geq 0.\label{kun}\]The above isomorphism is induced by the shuffle map and  following natural map \cite{L}
\[ \mathfrak{I}: [\alpha]\otimes[\beta]\longmapsto [\alpha \otimes \beta].\label{kun1}\] The same maps induce an  isomorphism for Hochschild homology of cyclic modules.

In the case of cyclic homology we do not have an isomorphism similar to \eqref{kun}. Instead, there is a short exact sequence \cite{bn, ka,kas1,L}:
$$\begin{CD}
  0  \longrightarrow Tot_{n}\mathcal{B}(C\otimes C') @>\mathcal{I} >> \displaystyle \bigoplus_{{i+j=n}}Tot_{i}\mathcal{B}(C)\otimes Tot_{j}\mathcal{B}(C')
   @>S\otimes id-id\otimes S >>  \end{CD}$$
   \[\begin{CD}
   \displaystyle \bigoplus_{{i+j=n-2}}Tot_{i}\mathcal{B}(C)\otimes Tot_{j}\mathcal{B}(C')\longrightarrow 0 .
   \end{CD}\label{kk}\] The map $\mathcal{I}$ is called the K\"{u}nneth map which can be defined as follows. Let $Tot\mathcal{B}(C)=k[u]\otimes C$, where $|u|=2 $, $Tot\mathcal{B}(C')=k[u']\otimes C'$ where $|u'|=2$,  and $Tot\mathcal{B}(C\otimes C')=k[v]\otimes C  \otimes C'$ where $|v|=2.$ One can define the map $\mathcal{I}$ as \cite{bn,L}:
   \begin{equation}\label{formula-kunneth}
     \mathcal{I}(v^{n})=\sum_{{p+q=n}}u^{p}u'^{q}.
     \end{equation}
      From the K\"{u}nneth short exact sequence \eqref{kk}, we get the following  long exact sequence in homology:

   $$\begin{CD}
    ... \longrightarrow HC_{n}(C\otimes C')@>\mathcal{I} >> \displaystyle \bigoplus_{{i+j=n}} HC_{i}(C)\otimes HC_{j}(C')
    @>S\otimes id-id\otimes S >>  \end{CD}$$
    $$\begin{CD}
    \displaystyle \bigoplus_{{i+j=n-2}}HC_{i}(C)\otimes HC_{j}(C')@> \partial >>
    HC_{n-1}(C\otimes C')\longrightarrow ....
    \end{CD}$$ Similarly for cyclic cohomology, we have the following short exact sequence:

$$\begin{CD}
0  \longrightarrow  Tot^{n}B(C \otimes C')  @> \mathcal{I} >> \displaystyle \bigoplus_{{i+j=n}}Tot^{i}B(C)\otimes Tot^{j}B(C') @> S\otimes id-id\otimes S >> \end{CD}$$

$$\begin{CD}
\displaystyle \bigoplus_{{i+j=n+2}}Tot^{i}B(C)\otimes Tot^{j}B(C')\longrightarrow 0. \end{CD}$$ Consequently, we get the following long exact sequence in cohomology \cite{bn, kas}:
$$ \begin{CD}
...\longrightarrow HC^{n}(C\otimes C') @> \mathcal{I} >> \displaystyle \bigoplus_{{p+q=n}}HC^{p}(C)\otimes HC^{q}(C')
@>S\otimes id-id \otimes S >>\end{CD}$$
 \begin{equation}
 \begin{CD}
 \displaystyle \bigoplus_{{r+s=n+2}}HC^{r}(C)\otimes HC^{s}(C')@>\partial >>HC^{n+2}(C\otimes C') \longrightarrow ....\end{CD}\end{equation}
where $\partial$ is the connecting morphism of short exact sequences.

 Now we state the K\"{u}nneth formula for the periodic case. One can find the formula for algebras in \cite{em}. Here we provide a similar argument in the general case of (co)cyclic modules. Recall that a supercomplex is a $\mathbb{Z}_{2}$-graded vector space $V=V_{0}\oplus V_{1}$ endowed with linear maps $\partial_{0}:V_{0}\longrightarrow V_{1}$ and  $\partial_{1}:V_{1}\longrightarrow V_{0}$  such that $\partial_{0}\partial_{1}=0$ and $\partial_{1}\partial_{0}=0$. We denote the corresponding chain complex by
$V_{0}   \overset{\partial_{0}}{\underset{\partial_{1}}{\leftrightarrows}} V_{1}$, where its homology is a $\mathbb{Z}_{2}$- graded vector space given by:
 $$ H_{0}=Ker\partial_{0}/Im\partial_{1},  \qquad       H_{1}=Ker\partial_{1}/Im\partial_{0}. $$ For example the inverse limit $ \underleftarrow{\lim} Tot \mathcal{B}(C) $ is a supercomplex. In this case $ V=TotB_{n}(C)$ , $\partial=B + b$ and we have $\underleftarrow{\lim}Tot\mathcal{B}(C)= (\prod C_{2n})_{n\geq 0}\oplus (\prod C_{2n+1})_{n\geq 0}$. For cohomology,  $ \underrightarrow{\lim}{Tot\mathcal{B}(C)}= V_{0}\oplus V_{1}=(\oplus C^{2n})_{n\geq 0} \oplus (\oplus C^{2n+1})_{n\geq 0}$. A map of supercomplexes  is a linear map $f:V=V_{0}\oplus V_{1} \longrightarrow W=W_{0}\oplus W_{1}$ such that sends $V_{0}$ to $W_{0}$ and  $V_{1}$ to $W_{1}$. If we have two supercomplexes $V$ and $W$, then $V\widehat{\otimes}W$ is a supercomplex where $(V\widehat{\otimes}W )_{0}=(V_{0}\otimes W_{0})\oplus (V_{1}\otimes W_{1})$ and $(V\widehat{\otimes}W )_{1}=(V_{0}\otimes W_{1})\oplus (V_{1}\otimes W_{0})$. We define

\begin{displaymath}
{\partial_{0}^{V\widehat{\otimes}W}} =
\left(\begin{array}{c c}
1\otimes \partial_{0}^{W} &  \partial_{1}^{V}\otimes 1\\

\partial_{0}^{V}\otimes 1 & -1\otimes \partial_{1}^{W}
\end{array}\right)
\qquad and \qquad
{\partial_{1}^{V\widehat{\otimes}W}} =
\left(\begin{array}{c c}
1\otimes \partial_{1}^{W} &  \partial_{1}^{V}\otimes 1\\

\partial_{0}^{V}\otimes 1 & -1\otimes \partial_{0}^{W}
\end{array}\right)
\end{displaymath}

The K\"{u}nneth formula for supercomplexes  holds:
 $$H(X \widehat{\otimes}Y)\cong H(X)\widehat{\otimes}H(Y).$$  We define the supercomplex map:

  \[\nabla:(\underleftarrow{\lim}Tot\mathcal{B}(\overline{C }))\widehat{\otimes} (\underleftarrow{\lim}Tot\mathcal{B}(\overline{C'}))\longrightarrow \underleftarrow{\lim}Tot\mathcal{B}(\overline{C}\otimes \overline{C'}).\] The restriction of $\nabla$ on $ (\prod_{n\geq 0} C^{2n})\otimes (\prod_{n\geq 0} C'^{2n})$ sends $\{\xi_{n}\}_{n}\otimes \{\xi'_{n}\}_{n} $ onto $\{\sum_{i+j=n}\xi_{i}\otimes \xi_{j}\}_{n}$. Similarly one can define the restriction of $I$ on $(\prod_{n\geq 0} C^{2n+1})\otimes (\prod_{n\geq 0} C'^{2n+1})$, $(\prod_{n\geq 0} C^{2n})\otimes (\prod_{n\geq 0} C'^{2n+1})$ and $(\prod_{n\geq 0} C^{2n+1})\otimes (\prod_{n\geq 0} C'^{2n})$.
  To prove that $\nabla$ is a quasi-isomorphism, we have two major problems. The first one is the fact that generally  homology does not commute with inverse limit, i.e.,
  $ H_{*}(\underleftarrow{\lim}Z_{m})\neq \underleftarrow{\lim}H_{*}(Z_{m}),$  where $Z_{m}$ is an inverse system.
   To solve this problem we need the Mittag-Leffler condition which guarantees the commutativity. The second problem is in general $(\prod_{{i=1}}^{{\infty}}V_{i})\otimes W\ncong \prod_{{i=1}}^{{\infty}}(V_{i}\otimes W),$  where $V_{i}$ and $W$ are some vector spaces . If $W$ is  finite dimensional then we obtain an isomorphism . In our case we can think about  $V_{i}$ and $W$ as the homology of a cyclic module. Now we are ready to have the following theorem:

   \bth \label{har}Suppose $C$ is a cyclic module which has the following two properties:

   (i) The inverse system $(HC(C)[-2m],S)_{m}$ satisfies the Mittag-Leffler condition.

   (ii) The periodic cyclic homology  $HP_{*}(C)$ is a finite dimensional vector space.

   Then, the map
   \[\nabla: HP_{*}(C)\widehat{\otimes}HP_{*}(C')\longrightarrow HP_{*}(C\otimes C')\]   is an isomorphism for any cyclic module $C'$.
   \ethe
   \bpf The proof is similar to the one in \cite{em} for algebras. Here we just mention that the condition (i) gives us $HP(C)= \underleftarrow{\lim}HC(C)[-2m]$ and (ii) proves $$(HP(C)\widehat{\otimes} \prod HC(C')[-2m])_{*} \cong \prod (HP_{*-[m]}(C)\otimes HC_{m}(C')).$$
    \epf

Thus under the assumptions of the Theorem \ref{har}, one obtains the K\"{u}nneth formula for periodic cyclic homology as follows:

   \[HP_{0}(C)\otimes HP_{0}(C') \oplus HP_{1}(C)\otimes HP_{1}(C') \cong HP_{0}(C\otimes C'),\label{per}\] and
   \[HP_{0}(C)\otimes HP_{1}(C') \oplus HP_{1}(C)\otimes HP_{0}(C') \cong HP_{1}(C\otimes C').\label{per1}\]

   For periodic cyclic cohomology, since direct limit commutes with cohomology, we do not need (i) and also since $\underrightarrow{\lim}Tot\mathcal{B}(C)=(\oplus C^{2n})_{n\geq 0}\oplus (\oplus C^{2n+1})_{n\geq 0})$ and $\oplus_{i\geq 0} (V_{i}\otimes W)\simeq (\oplus_{i\geq 0} V_{i})\otimes W$, for all vector spaces $V_{i}$ and $W$, we do not need (ii). Therefore we have the following theorem:

 \bth  Suppose $C$ and $C'$ are two  cocyclic modules.
   Then, the map
   \[\nabla: HP^{*}(C)\widehat{\otimes}HP^{*}(C')\longrightarrow HP^{*}(C\otimes C'),\label{arr}\]  is an isomorphism and we have the relations \eqref{per}  and \eqref{per1} in cohomology.

   \ethe

\section{Eilenberg-Zilber isomorphisms and K\"{u}nneth formulas in Hopf cyclic theory}
In this section we establish  the Eilenberg-Zilber isomorphisms and K\"{u}nneth formulas for Hopf cyclic theory with coefficients. Throughout the paper we assume that $\mathcal{H}$ is a Hopf algebra with a bijective antipode over a field $k$ of characteristic zero. The  coproduct, counit and antipode of $\mathcal{H}$ are denoted by
$\Delta$,  $\varepsilon$ and $S$,  respectively. For the coproduct we use the Sweedler notation in the form $\Delta(h)=h^{(1)}\ot h^{(2)}$; for a left coaction of $\mathcal{H}$ on a comodule $M$, $\blacktriangledown:M\longrightarrow \mathcal{H}\otimes M$, we write $\blacktriangledown(m)=m^{(-1)}\ot m^{(0)}$, and for a right coaction we write  $\blacktriangledown(m)=m^{(0)}\ot m^{(1)}$. If $f$ is a map of (co)chain complexes, then $f^{*}$ denotes the induced map on (co)homology.

\bde \cite{hkrs1} Let $\mathcal{H}$ be a Hopf algebra and $M$ a right module and a left comodule over $\mathcal{H}$. We call $M$ an anti-Yetter-Drinfeld module, if the action and coaction are compatible in the following sense:
$$(mh)^{(-1)}\otimes (mh)^{(0)} = S(h^{(3)})m^{(-1)}h^{(1)}\otimes m^{(0)}h^{(2)},$$
for all $m \in M$ and $h\in \mathcal{H}$. We call $M$ stable if for all $m\in M$,  $m^{(-1)}m^{(0)}=m $.
\ede

We abbreviate stable anti-Yetter-Drinfeld by SAYD. For example, given a character, i.e., a unital algebra map $\delta: \mathcal{H}\longrightarrow k$, and a group-like element
$\sigma \in \mathcal{H}$, one defines a SAYD module  $ M=^{\sigma}\!\!\!k_{\delta},$ where the action of $\mathcal{H}$ is defined by the character $\delta$, $mh=\delta(h)m,$ and the coaction via the group-like element $\sigma$, $\blacktriangledown(m)=\sigma\otimes m$. Then $ M=^{\sigma}\!\!\!k_{\delta}$ is  stable  if and only if $(\delta, \sigma)$ is a modular pair, i.e., $\delta(\sigma)=1$, and anti-Yetter-Drinfeld if and only if $(\delta, \sigma)$ is in involution, i.e.,
 $\sig^{-1} \ws^2(h) \sig = h,  \quad \forall h \in \mathcal{H}$, where $\ws(h):=\del(h^{(1)})S(h^{(2)})$\cite{cm3}. This module is called a modular pair in involution.

Given a Hopf algebra $\mathcal{H}$ equipped with a SAYD  $\mathcal{H}$-module $M$, we have the cocyclic module $C(\mathcal{H}, M)$ where $C^{n}(\mathcal{H}, M)= M\otimes_{\mathcal{H}} \mathcal{H}^{\otimes (n+1)},n\geq 0.$ Its cofaces, codegeneracies and cocyclic map are as follows \cite{hkrs2}:
\begin{eqnarray}
\delta_i(m\ot_{\mathcal{H}} h_0 \otimes \dots \otimes h_{n-1})
\!\!\!&=&\!\!\!
m\ot_{\mathcal{H}}  h_0 \otimes\dots
\otimes  h_i^{(1)}\otimes h_i^{(2)}\dots\otimes h_{n-1},
 ~~~ 0 \leq i <n ~~~\mbox{(cofaces)},~~~~~~~~~
\nonumber\\
\delta_{n}(m\ot_{\mathcal{H}} h_0 \otimes \dots \otimes h_{n-1})
\!\!\!&=&\!\!\!
m^{(0)}\ot_{\mathcal{H}}  h_0^{(2)}\otimes h_1
\otimes \dots \otimes h_{n-1}  \otimes m^{({-1})}h_0^{(1)}
~~~ \mbox{(flip-over face)},~~~~~~~~~
\nonumber\\
\sigma_i(m\ot{_\mathcal{H}} h_0 \otimes \dots \otimes h_{n+1})
\!\!\!&=&\!\!\!
m\ot_{\mathcal{H}} h_0 \otimes \dots
 \otimes \varepsilon(h_{i+1})\otimes
\dots\otimes h_{n+1},~~~0\leq i \leq n~\mbox{(codegeneracies)},~~~~~~~~~
\nonumber\\
\tau_n(m\ot_{\mathcal{H}} h_0 \otimes  \dots \otimes h_n)
\!\!\!&=&\!\!\!
m^{(0)}\ot_{\mathcal{H}} h_1 \otimes
\dots \otimes h_n \otimes m^{({-1})}h_0 ~~~ \mbox{(cocyclic map)}.~~~~~~~~~\nonumber
\end{eqnarray}

For $ M=^{\sigma}\!\!\!k_{\delta}$, the complex
$C(\mathcal{H},M)$ reduces to the cocyclic module of Connes-Moscovici \cite{cm3}.

\ble If $\mathcal{H}$ and $\mathcal{K}$ are Hopf algebras and $M$ and $N,$  SAYD modules over $\mathcal{H}$ and $\mathcal{K}$ respectively, then $ M\otimes N $ is a SAYD module over $ \mathcal{H}\otimes \mathcal{K}$ in a natural way.
\ele

\bpf We define the right action by:
\[ (m \otimes n)(h\otimes k)\ = mh\otimes nk,\] and the left coaction by:
\[ \blacktriangledown( m \otimes n ) = m^{(-1)}\otimes n^{(-1)}\otimes m^{(0)}\otimes n^{(0)}. \] We check the compatibility of action and coaction:
 \begin{eqnarray}
 &&\blacktriangledown(( m \otimes n )(h\otimes k)) = \blacktriangledown(  mh\otimes nk )= (mh\otimes nk)^{(-1)}\otimes (mh\otimes nk)^{(0)} \nonumber \\
 &=& (mh)^{(-1)}\otimes (nk)^{(-1)}\otimes (mh)^{(0)}\otimes (nk)^{(0)} \nonumber \\
 &=& S(h^{(3)})m^{(-1)}h^{(1)}\otimes S(k^{(3)})n^{(-1)}k^{(1)}\otimes m^{(0)}h^{(2)} \otimes n^{(0)}k^{(2)} \nonumber \\
 &=& S_{\mathcal{H}\otimes \mathcal{K}}(h^{(3)}\otimes k^{(3)})(m\otimes n)^{(-1)}(h^{(1)}\otimes (k^{(1)})\otimes ((m\otimes n)^{(0)} )(h^{(2)}\otimes k^{(2)}).\nonumber
\end{eqnarray}
 To check the stability:
 $$(m^{(0)}\otimes n^{(0)} )(m^{(-1)}\otimes n^{(-1)} ) = (m^{(0)}m^{(-1)}\ot n^{(0)}n^{(-1)}) = (m\ot n).$$
\epf

We need the following statement later.

 \ble If $\mathcal{H}$ and $\mathcal{K}$ are Hopf algebras and $M$ and $N$  SAYD modules over $\mathcal{H}$ and $\mathcal{K}$ respectively, then the following map
$$ \Omega_{r} : (M\otimes N)\otimes_{\mathcal{H}\otimes \mathcal{K}} (\mathcal{H}\otimes \mathcal{K})^{\otimes r+1}\longrightarrow (M\otimes_{\mathcal{H}} \mathcal{H}^{\otimes r+1})\otimes (N\otimes_{\mathcal{K}} \mathcal{K}^{\otimes r+1}),\quad r\geq 0,$$ given by $$ ((m\otimes n)\otimes_{{\mathcal{H}\ot \mathcal{K}}} (h_{0}\otimes k_{0})\otimes\dots\otimes(h_{r}\otimes k_{r}))\longmapsto (m\otimes{_\mathcal{H}} h_{0}\otimes\dots\otimes h_{r})\otimes(n\otimes_{\mathcal{K}} k_{0}\otimes \dots\otimes k_{r}),$$ is well-defined and defines an isomorphism of cocyclic modules: $$C(\mathcal{H}\otimes \mathcal{K}, M\otimes N)\cong  C(\mathcal{H}, M)\times C(\mathcal{K}, N).$$
\ele
\bpf Since we consider $\mathcal{H}$, $\mathcal{K}$ and $\mathcal{H}\ot \mathcal{K}$ as $\mathcal{H}^{\ot m}$, $\mathcal{K}^{\ot m}$ and $(\mathcal{H}\ot\mathcal{K})^{\ot m}$-modules, respectively,  by diagonal action, the map $\Omega$ is well-defined. We show that $\Omega_{r}$ commutes with $\delta_{i} $ and $\tau_{r},$ where $0\leq i < r $. The commutativity of $\Omega_{r}$ with $\delta_{r}$ and $\sigma_{i}$'s is left to the reader. For $\delta_{i}$'s we have:
 \begin{eqnarray}
 &&\Omega_{r} \delta_{i}(( m\otimes n)\otimes_{{\mathcal{H}\ot \mathcal{K}}} (h_{0}\otimes k_{0})\otimes \dots \otimes (h_{r-1}\otimes k_{r-1}))  \nonumber \\
 &=& \Omega_{r}((m\otimes n)\otimes_{{\mathcal{H}\ot \mathcal{K}}} (h_{0}\otimes k_{0})\otimes\dots\ot (h_{i}\otimes k_{i})^{(1)}\otimes (h_{i}\otimes k_{i})^{(2)}\otimes\dots\otimes (h_{r-1}\otimes k_{r-1})\nonumber \\
 &=& \Omega_{r}((m\otimes n)\otimes_{{\mathcal{H}\ot \mathcal{K}}} (h_{0}\otimes k_{0})\otimes\dots \ot(h_{i}^{(1)}\otimes k_{i}^{(1)})\otimes (h_{i}^{(2)}\otimes k_{i}^{(2)})\otimes\dots(h_{r-1}\otimes k_{r-1}))\nonumber \\
 &=& (m\otimes_{\mathcal{H}} h_{0}\otimes \dots \otimes h_{i}^{(1)}\otimes h_{i}^{(2)}\otimes\dots \ot h_{r-1})\otimes
 (n\otimes_{{\mathcal{K}}} k_{0}\otimes \dots \otimes k_{i}^{(1)}\otimes k_{i}^{(2)}\otimes\dots \ot k_{r-1}) \nonumber \\
 &=& (\delta_{i}\otimes \delta_{i}) (m\otimes_{{\mathcal{H}}} h_{0}\otimes \dots\otimes h_{r-1})\otimes (n\otimes_{{\mathcal{K}}} k_{0}\otimes \dots\otimes k_{r-1})\nonumber \\
&=&(\delta_{i}\otimes \delta_{i}) \Omega_{r}((m\otimes n)\otimes_{{\mathcal{H}\ot \mathcal{K}}} (h_{0}\otimes k_{0})\otimes \dots \otimes (h_{r-1}\otimes k_{r-1})). \nonumber
\end{eqnarray}
For $\tau_{r}$:
\begin{eqnarray}
 &&\Omega_{r} \tau_{r}((m\otimes n)\otimes_{{\mathcal{H}\ot \mathcal{K}}}  (h_{0}\otimes k_{0})\otimes \dots \otimes (h_{r}\otimes k_{r}))\nonumber \\
&=& \Omega_{r}( (m^{(0)}\otimes n^{(0)})\otimes_{{\mathcal{H}\ot \mathcal{K}}}  (h_{0}\otimes k_{0})\otimes\dots\otimes (h_{r}\otimes k_{r})\otimes ((m^{(-1)}\otimes n^{(-1)})(h_{0}\otimes k_{0})) \nonumber \\
&=&(m^{(0)}\otimes_{{\mathcal{H}}} h_{1}\otimes \dots\otimes h_{r}\otimes m^{(-1)}h_{0})\otimes (n^{(0)}\otimes_{{\mathcal{K}}} k_{1}\otimes \dots\otimes k_{r}\otimes n^{(-1)}k_{0})\nonumber \\
&=&(\tau_{r}\otimes \tau_{r})((m\otimes_{{\mathcal{H}}} h_{0}\otimes \dots\otimes h_{r})\otimes (n\otimes_{{\mathcal{K}}} k_{0}\otimes \dots\otimes k_{r})\nonumber \\
&=&(\tau_{r}\otimes \tau_{r})\Omega_{r}(((m\otimes n)\otimes_{{\mathcal{H}\ot \mathcal{K}}}  (h_{0}\otimes k_{0})\otimes \dots \otimes (h_{r}\otimes k_{r}))
\nonumber .
\end{eqnarray}
The bijectivity of $\Omega_{r}$ is obvious.
\epf

Using the Eilenberg-Zilber isomorphisms \eqref{EZIH}, \eqref{EZIC}, \eqref{abcde}, K\"{u}nneth formulas \eqref{kun},  \eqref{arr}, coupled with the above two lemmas, we obtain the following propositions.

\bpr Let $\mathcal{H}$ and $\mathcal{K}$ be two Hopf algebras and $M$ and $N$ be SAYD modules over $\mathcal{H}$ and $\mathcal{K}$ respectively. We have the the following  isomorphisms for Hopf Hochschild, cyclic, and periodic cyclic cohomology of Hopf algebras with coefficients:

\[HH^{*}(\mathcal{H}\otimes \mathcal{K}, M\otimes N)\cong HH^{*}(C(\mathcal{H},M)\otimes C(\mathcal{K},N)),\label{EHH}\]
\[HC^{*}(\mathcal{H}\otimes \mathcal{K}, M\otimes N)\cong HC^{*}(C(\mathcal{H},M)\otimes C(\mathcal{K},N)),\]
\[HP^{*}(\mathcal{H}\otimes \mathcal{K}, M\otimes N)\cong HP^{*}(C(\mathcal{H},M)\otimes C(\mathcal{K},N)).\]
\epr

\bpr
Let $\mathcal{H}$ and $\mathcal{K}$ be two Hopf algebras and $M$ and $N$ be SAYD modules over $\mathcal{H}$ and $\mathcal{K}$ respectively. We have the following isomorphisms of vector spaces,

\[ HH^{n}(\mathcal{H}\otimes \mathcal{K},M\otimes N)\cong \bigoplus_{{i+j=n}} HH^{i}(\mathcal{H},M)\otimes HH^{j}(\mathcal{K},N),\label{KHH}\]
and
\[HP^{n}(\mathcal{H}\otimes \mathcal{K},M\otimes N)\cong \bigoplus_{{i+j=n}} HP^{i}(\mathcal{H},M)\otimes HP^{j}(\mathcal{K},N),\quad n=0,1.\label{KHP}\]
and the following long exact sequence of differential graded vector spaces,
$$ \begin{CD}
...\longrightarrow HC^{n}(\mathcal{H}\otimes \mathcal{K},M\otimes N) @> \mathcal{I} >> \displaystyle \bigoplus_{{p+q=n}}HC^{p}(\mathcal{H},M)\otimes HC^{q}(\mathcal{K},N)
@>S\otimes id-id \otimes S >>\end{CD}$$
$$ \begin{CD}
 \displaystyle \bigoplus_{{r+s=n+2}}HC^{r}(\mathcal{H},M)\otimes HC^{s}(\mathcal{K},N)@>\partial >>HC^{n+2}(\mathcal{H}\otimes \mathcal{K},M\otimes N)\longrightarrow .... \end{CD}$$
\epr
\section{Coproducts in Hopf cyclic cohomology}
In this section we define  coproducts for Hopf Hochschild, cyclic and periodic cyclic cohomology of cocommutative Hopf algebras. The following statement plays an important role in the definition of these coproducts.

\bpr  Let $\mathcal{H}$ be a cocommutative Hopf algebra and  $M$  a SAYD module over $\mathcal{H}$. Let
$$\Phi_{n}= \psi\otimes \Delta^{\otimes n+1}: C^{n}(\mathcal{H},M)\longrightarrow C^{n}(\mathcal{H}\otimes \mathcal{H}, M\otimes M), $$ be a linear map, where $\psi: M\longrightarrow M\otimes M$, satisfying the following condition
\[(\Delta \otimes \psi) \circ \blacktriangledown_{M} = \blacktriangledown_{M\otimes M}\circ \psi. \label{shart}\] Then  the map
$$\rho_{n}=\Omega_{n} \Phi_{n} : C^{n}(\mathcal{H}, M)\longrightarrow (C(\mathcal{H},M)\times C(\mathcal{H},M))^{n},$$ is a map of cocyclic modules.
\epr
\bpf Since $\mathcal{H}$ is cocommutative,  one can easily see that $\rho_{n}\delta_{i}=(\delta_{i}\otimes \delta_{i})\rho_{n}$ where  $0\leq i \leq n-1$. We show that $\rho_{n}$ commutes with $\delta_{n}$ . For this we use the summation notation $\psi(m)= m_{(1)} \otimes m_{(2)}$. The condition \eqref{shart} is equivalent to:
$$ (m_{(1)})^{(-1)}\otimes (m_{(2)})^{(-1)} \otimes (m_{(1)})^{(0)}\otimes (m_{(2)})^{(0)}= m^{(-1)(1)}\otimes m^{(-1)(2)}\otimes (m^{(0)})_{(1)}\otimes (m^{(0)})_{(2)}.$$ By cocommutativity of $\mathcal{H}$ we have:
\begin{eqnarray}
&&(\delta_{n}\otimes \delta_{n}) \Omega_{n-1} (\psi \otimes \Delta^{\otimes n})(m \otimes h_{0}\otimes \dots \otimes h_{n-1})\nonumber \\
&=&(\delta_{n}\otimes \delta_{n})\Omega_{n-1}(\psi(m)\otimes \Delta(h_{0})\otimes \dots \otimes \Delta(h_{n-1}))\nonumber \\
&=&(\delta_{n}\otimes \delta_{n})(m_{(1)}\otimes h_{1}^{(1)}\otimes \dots h_{{n-1}}^{(1)})\otimes (m_{(2)}\otimes  h_{1}^{(2)}\otimes \dots h_{{n-1}}^{(2)})\nonumber \\
&=&\Omega_{n-1}((m_{(1)})^{(0)}\otimes h_{0}^{(1)(2)}\otimes h_{2}^{(1)}\otimes \dots\otimes h_{n-1}^{(1)}\otimes  (m_{(1)})^{(-1)}h_{0}^{(1)(1)})\otimes ((m_{(2)})^{(0)}\otimes h_{0}^{(2)(2)}\otimes \nonumber \\
&&h_{1}^{(2)}\otimes \dots\otimes h_{n-1}^{(2)}\otimes (m_{(2)})^{(-1)}h_{0}^{(2)(1)})\nonumber \\
&=&\Omega_{n-1}(((m_{(1)})^{(0)}\otimes (m_{(2)})^{(0)})\otimes (h_{0}^{(1)(2)}\otimes h_{0}^{(2)(2)})\otimes (h_{1}^{(1)}\otimes h_{1}^{(2)})\otimes\dots\otimes \nonumber \\
&&(h_{n-1}^{(1)}\otimes h_{n-1}^{(2)})\otimes (m_{(1)})^{(-1)}h_{0}^{(1)(1)}\otimes (m_{(2)})^{(-1)}h_{0}^{(2)(1)}))\nonumber \\
&=&\Omega_{n-1}(\psi(m)^{(0)}\otimes (\Delta(h_{1}))^{(2)}\otimes \Delta(h_{2})\otimes\dots\otimes \Delta(h_{{n-1}})\otimes \psi(m)^{(-1)}(\Delta(h_{1}))^{(1)})\nonumber \\
&=&\Omega_{n-1}(\psi(m^{(0)})\otimes \Delta(h_{0}^{(2)})\otimes \Delta(h_{1})\otimes\dots \otimes \Delta(h_{n-1})\otimes \Delta(m^{(-1)})\Delta(h_{0}^{(1)}))\nonumber \\
&=&\Omega_{n-1}(\psi\otimes\Delta^{\otimes n})(m^{(0)}\otimes h_{0}^{(2)}\otimes h_{1}\otimes\dots\otimes h_{n-1}\otimes m^{(-1)}h_{0}^{(1)})\nonumber \\
&=&\Omega_{n-1}(\psi\otimes\Delta^{\otimes n})\delta_{n}(m\otimes h_{0}\otimes \dots\otimes h_{n-1}).\nonumber
\end{eqnarray}
We show that $\rho_{n}$ commutes with $\tau_n$:
\begin{eqnarray}
&&\Omega_{n}(\psi\otimes \Delta^{\otimes n+1 })\tau_n (m\otimes h_0\otimes h_1 \dots\otimes h_n )\nonumber \\
&=&\Omega_{n}(\psi\otimes\Delta^{\otimes n+1})(m^{(0)} \otimes h_1 \otimes \dots \otimes h_n \otimes m^{(-1)}h_0 ) \nonumber \\
&=&\Omega_{n}(\psi(m^{0})) \otimes \Delta(h_1) \otimes\dots \otimes \Delta(h_n)\otimes \Delta(m^{(-1)}\Delta(h_0))\nonumber \\
&=&\Omega_{n}(((m^{(0)})_{(1)}\otimes (m^{(0)})_{(2)})\otimes (h_{1}^{(1)}\otimes h_{1}^{(2)})\otimes \dots\otimes (h_{n}^{(1)}\otimes h_{n}^{(2)})\otimes (m^{(-1)(1)}h_{0}^{(1)}\otimes m^{(-1)(2)}h_{0}^{(2)}))\nonumber \\
&=&((m^{(0)})_{(1)}\otimes h_{1}^{(1)}\otimes\dots\otimes h_{n}^{(1)}\otimes m^{(-1)(1)}h_{0}^{(1)})\otimes ((m^{(0)})_{(2)}\otimes h_{1}^{(2)}\otimes\dots\otimes h_{n}^{(2)}\otimes m^{(-1)(2)}h_{0}^{(2)})\nonumber \\
&=&((m_{(1)})^{(0)}\otimes h_{1}^{(1)}\otimes\dots\otimes h_{n}^{(1)}\otimes m_{(1)}^{(-1)}h_{0}^{(1)})\otimes ((m_{(2)})^{(0)}\otimes h_{1}^{(2)}\otimes\dots\otimes h_{n}^{(2)}\otimes m_{(2)}^{(-1)}h_{0}^{(2)})\nonumber \\
&=&\tau_{n}(m_{(1)}\otimes h_{0}^{(1)}\otimes \dots\otimes h_{n}^{(1)})\otimes \tau_{n}(m_{(2)}\otimes h_{0}^{(2)}\otimes \dots\otimes h_{n}^{(2)})\nonumber \\
&=& (\tau_{n}\otimes \tau_{n})\Omega_{n}(\psi(m)\otimes \Delta(h_{0})\otimes \dots\otimes \Delta(h_{n}))\nonumber\\
&=&(\tau_n \ot \tau_n)\Omega_n (\psi\ot \Delta^{\ot(n+1)})(m\ot_{\mathcal{H}}h_{0}\ot\cdots \ot h_n)).
\end{eqnarray}
So we have:
$$ \Omega_{n}( \psi\otimes\Delta^{\otimes(n+1)})\tau_{n} = (\tau_{n}\otimes \tau_{n})\Omega_{n} (\psi\otimes \Delta^{\otimes(n+1)}).$$
\epf

Based on  the proof of the  previous lemma, in fact the map $\Phi_n$ is a map of cocyclic modules.
Now we have all the needed tools to define the desired coproducts.
The following theorem provides a coproduct for Hochschild cohomology of a cocommutative Hopf algebra with coefficients in a SAYD module.

\bpr Suppose $\mathcal{H}$ is a cocommutative Hopf algebra, $M$ a SAYD module over $\mathcal{H}$, and $\psi:M\longrightarrow M\otimes M, $  a linear map satisfying \eqref{shart}. The following map
\[\sqcup= (Sh \rho)^{*}: HH^{n}(\mathcal{H},M)\longrightarrow \displaystyle \bigoplus_{{i+j=n}}HH^{i}(\mathcal{H},M)\otimes HH^{j}(\mathcal{H},M),\label{grad}\]
defines a coproduct for Hopf Hochschild cohomology where $\rho$ is defined in Proposition $5.1$.
\epr
\bpf Since $H$ is cocommutative, by Proposition $5.1$ the map
$$\rho: C^n(H,M)\longrightarrow C^n(H,M)\ot C^n(H,M),$$ is a map of Hochschild complexes.
Also since the shuffle map $$Sh: C^n(H,M)\ot C^n(H,M)\longrightarrow (C(H,M)\ot C(H,M))^n, $$ is a  map of Hochschild complexes, the composition    $Sh \rho$ is a map of Hochschild complexes and therefore it induces the map  $(Sh \rho)^* $ on the level of cohomology. One notes that in this case the shuffle map $Sh$ is giving us the K\"{u}nneth morphism.
\epf

\bpr\label{prop} In the Connes-Moscovici case with $M= ^1\!\!\!k_\delta,$ we have an explicit formula for $Sh \rho$:

\[Sh_{n} \rho_{n}(h_{1}\otimes ...\otimes h_{n})=\sum_{{\sigma}}sign(\sigma) (h_{\sigma(1)}\otimes ...\otimes h_{\sigma(p)})\otimes (h_{\sigma(p+1)}\otimes ...\otimes h_{\sigma(p+q)}),\]
where $\sigma$ runs over all $(p,q)$-shuffles.
\epr
\bpf Using $\epsilon(h^{(1)})h^{(2)}=\epsilon(h^{(2)})h^{(1)}=h,$ we have:
\begin{eqnarray}
&&Sh_{n}\Omega_{n}\Phi_{n}(h_{1}\otimes...\otimes h_{n})\nonumber \\
&=& Sh_{n}((h_{1}^{(1)}\otimes...\otimes h_{n}^{(1)})\otimes (h_{1}^{(2)}\otimes...\otimes h_{n}^{(2)}))\nonumber \\
&=&\sum_{{\sigma}}sign(\sigma)\epsilon(h_{\sigma(1)}^{(2)})...\epsilon(h_{\sigma(p)}^{(2)})\epsilon(h_{\sigma(p+1)}^{(1)})...\epsilon(h_{\sigma(n)}^{(1)})(h_{\sigma(1)}^{(1)}\otimes...h_{\sigma(p)}^{(1)})\otimes (h_{\sigma(p+1)}^{(2)}\otimes...h_{\sigma(n)}^{(2)})\nonumber \\
&=&\sum_{{\sigma}} sign(\sigma)(\epsilon(h_{\sigma(1)}^{(2)})h_{\sigma(1)}^{(1)}\otimes...\otimes \epsilon(h_{\sigma(p)}^{(2)})h_{\sigma(p)}^{(1)})\otimes (\epsilon(h_{\sigma(p+1)}^{(1)})h_{\sigma(p+1)}^{(2)}\otimes...\otimes \epsilon(h_{\sigma(n)}^{(1)})h_{\sigma(n)}^{(2)})\nonumber \\
&=&\sum_{{\sigma}}sign(\sigma)(h_{\sigma(1)}\otimes...\otimes h_{\sigma(p)})\otimes (h_{\sigma(p+1)}\otimes...\otimes h_{\sigma(n)}).\nonumber
\end{eqnarray}
\epf

One knows that Hochschild homology of a commutative algebra is a graded commutative and associative algebra \cite{L}. Analogous to this classic result we have:

\bpr The Hopf Hochschild cohomology of a cocommutative Hopf algebra $\mathcal{H}$ with coefficients in a SAYD module $M$ equipped with map $\psi$  satisfying \eqref{shart}, is a graded cocommutative and coassociative coalgebra by \eqref{grad}.

\epr
\bpf  The cocommutativity and coassociativity can be verified by series of long, but straightforward, computation which we omit here.
\epf

\bth Suppose $\mathcal{H}$ is a cocommutative Hopf algebra and $M$ a SAYD module over $\mathcal{H},$ equipped with a map $\psi:M\longrightarrow M\otimes M$  satisfying \eqref{shart}. The following map
\[\sqcup= \mathcal{I}(\widetilde{Sh} \varrho)^{*}: HC^{n}(\mathcal{H},M)\longrightarrow \displaystyle \bigoplus_{p+q=n}HC^{p}(\mathcal{H},M)\otimes HC^{q}(\mathcal{H},M),\label{grad1}\]
defines a coproduct for Hopf cyclic cohomology of $\mathcal{H}$ with coefficients in $M$, where $\varrho_n= \oplus_{i\geq0} (\Omega_{n-2i}\Phi_{n-2i}),$ and $\mathcal{I}$ is defined in \eqref{formula-kunneth}.
\ethe
\bpf
The map $(\widetilde{Sh} \varrho)^{*}$ is induced by  composition of the following maps of complexes,

$$\begin{CD}Tot^nB(C(H,M))   @>\varrho_n >> Tot^nB(C(H,M)\times C(H,M)) @>{\widetilde{Sh_n}}>>  Tot^n B (C(H,M)\ot C(H,M))  \end{CD}$$

Since $H$ is cocommutative the maps $\Phi_{n-2i}$ are morphisms of cocyclic modules by Proposition $5.1$. Since $\Omega_{n-2i}$ are maps of cocyclic modules by Lemma $4.2$, the morphism   $\varrho_n$ is a map of cocyclic modules. Therefore the composition $\widetilde{Sh} \varrho$ is a map of cocyclic modules which induces the map $(\widetilde{Sh} \varrho )^*$ on the level of cyclic cohomology. Now we compose this map by K\"{u}nneth map $\mathcal{I}$  defined in \eqref{formula-kunneth} to obtain  the coproduct $\sqcup$ on the level of cyclic cohomology.

\epf

\bth Suppose $\mathcal{H}$ is a cocommutative Hopf algebra and $M $ a SAYD module over $\mathcal{H},$  equipped with a map $\psi:M\longrightarrow M\otimes M$   satisfying \eqref{shart}. The following map
\[\sqcup= \nabla(\widetilde{Sh} \varrho)^{*}: HP^{n}(\mathcal{H}, M)\longrightarrow  \bigoplus_{{i+j= n}} HP^{i}(\mathcal{H}, M)\otimes HP^{j}(\mathcal{H}, M),\quad n=0,1,\label{grad2}\]
defines a coproduct for periodic Hopf cyclic cohomology where $\overline{\varrho}=\oplus_{i\geq0}\Phi_{2i+ \ast},$ and $\nabla$ is as \eqref{KHP}.
\ethe

Now we provide some examples of the map $\psi$ satisfying the condition \eqref{shart}.

\bex Let $\mathcal{H}=kG $ be the group algebra of the discrete group $G$. Suppose $M$ is a SAYD module over $\mathcal{H}.$ One can check that $M$ is a $G$-graded vector space $M = \oplus_{g\epsilon G }M_{g} $, where the coaction $\blacktriangledown$ is defined by $\blacktriangledown(m)= g \otimes m. $ The stability condition implies $gm=m$ for all $ m \in M_{g},$  and anti-Yetter-Drinfeld condition is equivalent to  $hm \in M_{hgh^{-1}}$ for all $g$ and $h \in G$. Since $g$ is a group-like element, any linear map $\psi: M\longrightarrow M\otimes M,$ with $\psi(M_{g})\subseteq M_{g}\otimes M_{g}$ satisfies the condition \eqref{shart}.
\eex

\bex Let $\mathcal{H}=\mathcal{U}(\mathfrak{g})$ be the universal enveloping  algebra  of a Lie algebra $\mathfrak{g}$ and $M$ an arbitrary module over $\mathcal{H}$. We can define a comodule structure on $M$ by trivial coaction: $m \longmapsto 1\otimes m $. It can be shown that $M$ is a SAYD module over $\mathcal{H}$ and any linear map $\psi: M\longrightarrow M\otimes M, $ satisfies \eqref{shart}.
\eex

\bex  Let $\mathcal{H}$ be any cocommutative Hopf algebra and $M= ^1\!\!\!k_\delta$. It is easy to check that any linear map $\psi: M\longrightarrow M\otimes M, $ satisfies the condition \eqref{shart}. We use this example to get a coproduct in Connes-Moscovici setting in the following corollary.
\eex

\bco  If $M= ^1\!\!\!k_\delta$ then the coproducts  \eqref{grad}, \eqref{grad1}, \eqref{grad2} reduce to coproducts for Connes-Moscovici Hopf Hochschild, cyclic and periodic cyclic cohomology of a cocommutative Hopf algebra $\mathcal{H}$.
\eco

\section{Relation with coproduct in Lie algebra homology}

In this section we show that the coproduct  \eqref{grad2}  for periodic Hopf cyclic cohomology for $\mathcal{H}=\mathcal{U}(\mathfrak{g})$,  the universal enveloping algebra of a Lie algebra $\mathfrak{g}$,  agrees with the coproduct  in Lie algebra homology.
Recall that the periodic Hopf cyclic cohomology of $\mathcal{H}=\mathcal{U}(\mathfrak{g})$ is given by \cite{cm2}

\[HP_{(\delta,1)}^{\ast}(\mathcal{U}(\mathfrak{g}))\cong \bigoplus_{{k\geq 0}}H_{2k+\ast}(\mathfrak{g},k_{\delta}),\label{best}\]
 where $\ast=0,1$, and $k_{\delta}$ is a $\mathfrak{g}$-module via the character $\delta$. The right hand side of this isomorphism is the Lie algebra homology of $\mathfrak{g}$ with coefficients in $k_{\delta}$. It is the homology of the the following mixed complex:
 \[
\xymatrix{ \bigwedge^0\mathfrak{g}   \overset{d_{Lie}}{\underset{0}{\leftrightarrows}}  \bigwedge^1\mathfrak{g}   \overset{d_{Lie}}{\underset{0}{\leftrightarrows}}\bigwedge^2\mathfrak{g}   \overset{d_{Lie}}{\underset{0}{\leftrightarrows}}\dots }
,\] where $d_{Lie}$ denotes the Chevalley-Eilenberg differential for Lie algebra homology. The isomorphism \eqref{best} is induced by the anti symmetrization map $A_{n}:  \bigwedge^{n}\mathfrak{g} \longrightarrow \mathcal{U}(\mathfrak{g})^{\otimes n },$ given by $$ (g_{1}\wedge \dots\wedge g_{n}) \longmapsto \frac{1}{n!}(\sum_{\sigma}sign(\sigma)( g_{\sigma(1)}\otimes \dots\otimes    g_{\sigma(n)}).$$ Here $\sigma $ runs over all permutations of the set $\{1,2,\dots,n\}.$ If $\delta=\varepsilon,$ then $d_{Lie}$ is given by
$$ d_{Lie}(g_{1}\wedge\dots\wedge g_{n})= \sum_{1\leq i\leq i\leq n}(-1)^{i+j+1}[g_{i},g_{j}]\wedge (g_{1}\wedge\dots\wedge\hat{g_{i}}\wedge\dots\wedge\hat{g_{j}}\wedge\dots\wedge g_{n}). $$ Let $C_{n}^{Lie}(\mathfrak{g})= \bigwedge^{n}\mathfrak{g}$. One knows that Lie algebra homology with trivial coefficients is a coalgebra by the following coproduct

$$\cup_{Lie} :C_{n}^{Lie}(\mathfrak{g}) \longrightarrow \bigoplus_{{p+q=n}}(C_{p}^{Lie}(\mathfrak{g}))\otimes (C_{q}^{Lie}(\mathfrak{g})),$$
given by
\[  \cup_{Lie}(g_{1}\wedge\dots\wedge g_{n}) = \sum_{{\sigma}} sign(\sigma)(g_{\sigma(1)}\wedge \dots\wedge g_{\sigma(p)}) \otimes (g_{\sigma(p+1)}\wedge ...\wedge g_{\sigma(n)}),\] where $\sigma$ runs over  all $(p,q)$-shuffles. Since $ \varepsilon(g_{i})= 0 $ for all $ g_{i} \in  \mathfrak{g}$, the image of the anti symmetrization map is in the normalized complex of $C^{n}(\mathcal{U}(\mathfrak{g}))= \mathcal{U}(\mathfrak{g})^{\otimes n}$.
One can define the following map of mixed complexes:
\[\mathbb{A}_{n}= \sum_{i+j=n}A_{i}\otimes A_{j} : \bigoplus_{{i+j=n}} ( \Lambda^{i}\mathfrak{g})\otimes (\Lambda^{j}\mathfrak{g})\longrightarrow\bigoplus_ {{i+j=n}} (\mathcal{U}^{i}(\mathfrak{g}))\otimes (\mathcal{U}^{j}(\mathfrak{g})) .\]  It induces a map:
$$ \mathbb{A}^{*}:\bigoplus_{{k\geq0}} H_{2k+ \ast}(C^{Lie}(\mathfrak{g}))\otimes C^{Lie}(\mathfrak{g})))\longrightarrow HP^{\ast}(C(\mathcal{U}(\mathfrak{g})\otimes C(\mathcal{U}(\mathfrak{g})).$$

Now using the K\"{u}nneth formula \eqref{KHP}, we obtain the following map:

$$ \nabla(\mathbb{A}\cup_{Lie})^{*} : \bigoplus_{{n\geq0}}H_{2n+1}(\mathfrak{g},k_{\varepsilon})\longrightarrow HP^{1}_{(\varepsilon,1)}(\mathcal{U}(\mathfrak{g}))\otimes HP^{0}_{(\varepsilon,1)}(\mathcal{U}(\mathfrak{g})) \oplus HP^{0}_{(\varepsilon,1)}(\mathcal{U}(\mathfrak{g}))\otimes HP^{1}_{(\varepsilon,1)}(\mathcal{U}(\mathfrak{g})),$$ and similarly for the even case.

\bth Under the isomorphism \eqref{best}, the coproduct \eqref{grad2} for the periodic cyclic cohomology of the universal enveloping algebra $ \mathcal{U}(\mathfrak{g})$ with trivial coefficients coincides with the coproduct of Lie algebra homology. Equivalently, the following diagram commutes on the level of homology.

\[\begin{CD}
\bigoplus_{i=0} C_{2i+*}^{Lie}(\mathfrak{g}) @>A>> \bigoplus_{i=0} C^{2i+*}(\mathcal{U}(\mathfrak{g}))\\
@VV\cup_{Lie} V @VV\cup' V\\
\bigoplus_{i=0}( C^{Lie}(\mathfrak{g})\otimes C^{Lie}(\mathfrak{g}))_{2i+*}@>\mathbb{A}>> \bigoplus_{i=0} (C(\mathcal{U}(\mathfrak{g}))\otimes C(\mathcal{U}(\mathfrak{g})))^{2i+*}
\end{CD}\] \\

where $\cup'= \widetilde{Sh}\Omega\Delta^{\otimes n}$.

\ethe
\bpf The commutativity of the diagram is equivalent to
\[ sh_{n} \Omega_{n}  \Delta^{\otimes n}_{\mathcal{U}(\mathfrak{g})} A_{n} \oplus sh'_{n+2}\Omega_{n+2} \Delta^{\otimes n+2}_{\mathcal{U}(\mathfrak{g})} A_{n+2} = \mathbb{A}_{n}\cup_{Lie}, \quad n\in \mathbb{N}.\] Using the Proposition \ref{prop}, one can easily see
 $$sh_{n}\Omega_{n} \Delta^{\otimes n}_{\mathcal{U}(\mathfrak{g})}A_{n}= \mathbb{A}_{n}\cup_{Lie},$$ where

$$ \mathbb{A}_{n}\cup_{Lie}(g_{1} \wedge \dots \wedge g_{n}) = \sum_{p=0}^{n}\sum_{{\sigma \in S_{n}}} sign(\sigma)(g_{\sigma(1)}\otimes \dots\otimes g_{\sigma(p)})\otimes (g_{\sigma(p+1)}\otimes \dots\otimes g_{\sigma(n)}).$$ Now it is enough to show:

 \begin{align}
  sh'_{n+2}\Omega_{n+2} \Delta^{\otimes n+2}_{U(\mathfrak{g})} A_{n+2} = A_{n}d_{Lie},
  \end{align}

  which means this term will be zero on the level of homology.
 To do this, first we compute $\Omega_{n+2} \Delta_{\mathcal{U}(\mathfrak{g})}^{\otimes (n+2)}A_{n+2}:$

 \begin{eqnarray}
&&\Omega_{n+2} (\Delta^{\otimes (n+2)}_{\mathcal{U}( \mathfrak{g})})A_{n+2}(g_{1}\wedge \dots \wedge g_{n+2})\nonumber \\
&=& \Omega_{n+2} \sum_{{\sigma \in S_{n+2}}} \Delta^{\otimes (n+2) }_{\mathcal{U}(\mathfrak{g})}( g_{\sigma(1)}\otimes\dots \otimes g_{\sigma(n+2)}) \nonumber \\
&=& \Omega_{n+2} \sum_{{\sigma \in S_{n+2}}}(\Delta (g_{\sigma(1)})\Delta (g_{\sigma(2)})\dots \Delta (g_{\sigma(n+2)})) \nonumber \\
&=& \Omega_{n+2} \sum_{{\sigma \in S_{n+2}}}((1\otimes g_{\sigma(1)} + g_{\sigma(1)}\otimes 1)\dots(1\otimes g_{\sigma(n+2)} +g_{\sigma(n+2)}\otimes 1)) \nonumber \\
&=& \sum_{p=0}^{n+2}\sum_{{\sigma \in S_{n+2}}}( g_{\sigma(1)} \otimes\dots \otimes g_{\sigma(p)}\otimes 1\otimes\dots\otimes1)\otimes (g_{\sigma(p+1)}\otimes\dots \otimes  g_{\sigma(n+2)}\otimes 1\otimes\dots\otimes1),\nonumber
\end{eqnarray}
where  in the last sum, $g_{\sigma(1)},\dots,g_{\sigma(p)} $ and $g_{\sigma(p+1)},\dots,g_{\sigma(p+q)} $ appear in $(p,q)$-shuffle spots. For example when $n=1$, we have: \begin{eqnarray}
&&\Omega_{3}  \Delta^{\otimes 3}_{\mathcal{U}(\mathfrak{g})}A_{3}(g_{1}\wedge g_{2} \wedge g_{3}) =
(1\otimes 1\otimes 1)\otimes(g_{\sigma(1)}\otimes g_{\sigma(2)}\otimes g_{\sigma(3)}) \nonumber \\
&+&(1\otimes 1\otimes g_{\sigma(3)})\otimes (g_{\sigma(1)}\otimes g_{\sigma(2)}\otimes 1)+(1\otimes g_{\sigma(2)}\otimes 1)\otimes (g_{\sigma(1)}\otimes 1\otimes g_{\sigma(3)})\nonumber \\
&+& (1\otimes g_{\sigma(2)}\otimes g_{\sigma(3)})\otimes (g_{\sigma(1)}\otimes1 \otimes 1)+(g_{\sigma(1)}\otimes g_{\sigma(2)}\otimes 1)\otimes (1\otimes 1\otimes g_{\sigma(3)})\nonumber \\
&+& (g_{\sigma(1)}\otimes g_{\sigma(2)}\otimes g_{\sigma(3)})\otimes (1\otimes 1\otimes 1)+(g_{\sigma(1)}\otimes 1\otimes 1)\otimes (1\otimes g_{\sigma(2)}\otimes g_{\sigma(3)})\nonumber \\
&+& (g_{\sigma(1)}\otimes 1\otimes g_{\sigma(3)})\otimes (1\otimes g_{\sigma(2)}\otimes 1).\nonumber
\end{eqnarray}
Since $ sh'_{n} = \oplus_{{i+j = n+2}} sh'_{i,j}$, it is enough to show $ sh'_{i,j} = 0 ,$ for all $1\leq i,j\leq n+1$. For $ sh'_{i,j}$, we apply degeneracies, i.e.,  counit $\varepsilon$, $n+2$ times on the following elements:
 $$( g_{\sigma(1)} \otimes\dots \otimes g_{\sigma(p)}\otimes 1\otimes\dots\otimes1)\otimes (g_{\sigma(p+1)}\otimes\dots \otimes  g_{\sigma(n+2)}\otimes 1\otimes\dots\otimes1),$$ and then we apply $(s_{\sigma'(1)}\otimes \dots s_{\sigma'(q))}\otimes (s_{\sigma'(q+1)}\dots s_{\sigma'(n)})$. Now for the terms which are zero, there is nothing remained to prove.  Those which are not zero,  should be in the forms of
 $$ (g_{\sigma(1)}\otimes \dots\otimes g_{\sigma(i)})\otimes (g_{\sigma(i+1)}\otimes \dots\otimes g_{\sigma(n+2)}),$$ for some $1\leq i\leq n+2$.
 Here $\sigma$ runs over all permutations in $S_{n+2}$, the symmetric group with $n+2$ terms. Now we compute
 $$ \sigma_{i-p-1}^{i-1}\tau ^{p+1}( g_{\sigma(1)}\otimes \dots\otimes g_{\sigma(i)})\otimes \sigma_{j-q-1}^{j-1}\tau ^{q+1}(g_{\sigma(i+1)}\otimes \dots\otimes g_{\sigma(n+2)}).$$  Let $\sigma(i+1),\dots,\sigma(n)$ be fixed. For any $0 \leq k\leq n+2$, we have:
\begin{eqnarray}
&&\sigma^{i-1}_{i-k}\tau^{k}(g_{\sigma(1)}\otimes\dots\otimes g_{\sigma(i)})\nonumber \\
&=&  \sigma^{i-1}_{i-k}(\sum_{{\sigma}} S(g^{(i)}_{\sigma{(k)}})g_{\sigma(k+1)}\otimes S(g^{(i-1)}_{\sigma{(k)}})g_{\sigma(k+2)}\otimes \dots\otimes S(g^{(k)}_{\sigma{(k)}})\sigma \nonumber \\
&& \otimes S(g^{(k-1)}_{\sigma{(k)}})g_{\sigma(1)}\otimes\dots\otimes S(g^{(2)}_{\sigma(k)})g_{\sigma(2)}
\otimes \widetilde{S}(g^{(1)}_{\sigma{(k)}})\sigma g_{\sigma(i-1)})\nonumber \\
&=&  \sigma^{i-1}_{i-k}(\sum_{{\sigma}\in S_{p}}\sum_{{r=k+1}}^{{i-1}} g_{\sigma(k+1)}\otimes\dots g_{\sigma(r-1)}\otimes \nonumber \\
&& g_{\sigma(k)}g_{\sigma(r)}\otimes g_{\sigma(r+1)}\otimes\dots\otimes g_{\sigma(i)}\otimes 1\otimes g_{\sigma(1)}\otimes\dots\otimes g_{\sigma(i-1)} \nonumber \\
&+&   g_{\sigma(k+1)}\otimes \dots\otimes g_{\sigma(i)}\otimes g_{\sigma(k)}\otimes g_{\sigma(1)}\otimes\dots\otimes g_{\sigma(i-1)})\nonumber \\
&=&   \sum_{{\sigma \in S_{i}}}\sum_{{r=k+1}}^{p-1} g_{\sigma(k+1)}\otimes\dots\otimes g_{\sigma(r-1)}\otimes\nonumber \\
&& g_{\sigma(k)}g_{\sigma(r)}\otimes g_{\sigma(r+1)}\otimes\dots\otimes g_{\sigma(i)}\otimes g_{\sigma(1)}\otimes\dots\otimes g_{\sigma(i-1)}\nonumber \\
&=&   \sum_{{\sigma \in S_{i}},\sigma(v)>\sigma(w)} g_{\sigma(1)}\otimes \dots\otimes g_{\sigma(v)}g_{\sigma(w)}\otimes\dots\otimes g_{\sigma(i)} \nonumber \\
&+&    \sum_{{\sigma \in S_{i}},\sigma(v)>\sigma(w)} g_{\sigma(1)}\otimes \dots\otimes g_{\sigma(w)}g_{\sigma(v)}\otimes\dots\otimes g_{\sigma(i)} \nonumber \\
&=&   \sum_{{\sigma \in S_{i}},\sigma(v)>\sigma(w)} g_{\sigma(1)}\otimes \dots\otimes g_{\sigma(v)}g_{\sigma(w)}- g_{\sigma(w)}g_{\sigma(v)} \otimes\dots\otimes g_{\sigma(i)}\nonumber \\
&=&   \sum_{{\sigma \in S_{i}},\sigma(v)>\sigma(w)} g_{\sigma(1)}\otimes \dots\otimes [g_{\sigma(v)},g_{\sigma(w)}] \otimes\dots\otimes g_{\sigma(i)}\nonumber \\
&=&   A([g_{\sigma(v)},g_{\sigma(w)}]\wedge  g_{\sigma(1)}\wedge \dots\wedge  \hat{g}_{\sigma(w)}\wedge\dots\wedge  \hat{g  }_{\sigma(v)}\wedge\dots\wedge  g_{\sigma(i)})\nonumber \\
&=&  Ad^{Lie}( g_{1}\wedge\dots\wedge g_{i}).\nonumber
\end{eqnarray}
\epf

\section{Coproducts in   Hopf cyclic homology}
In this section we define coproducts for Hopf cyclic homology in the sense of  \cite{hkrs1, kr1}. We show that for group algebras our coproduct coincides with the coproduct in group homology. In this section  we assume $M$ is a left-left stable anti-Yetter-Drinfeld module on a Hopf algebra $\mathcal{H}$. Let
$$\widetilde{C}_{n}(\mathcal{H},M)= \mathcal{H}^{\otimes n+1}\Box_{\mathcal{H}}M \quad n\geq0,$$ where  $X\square_{\mathcal{H}}Y =\ker(\Delta_{X}\otimes id- id \otimes\,_Y\Delta): X\otimes Y\longrightarrow X\otimes \mathcal{H}\otimes Y$ is the cotensor product of $X$ and $Y$ \cite{hkrs1, kr1}.
One can define faces, degeneracies and cyclic maps on
$\{\widetilde{C}_{n}(\mathcal{H},M)\}_{n\in \mathbb{N}}$ as follows:
\bea
&&\
\delta_{i}(h_0\ot\dots\ot h_n\ot m)= h_0\ot \dots
 \ot h_i h_{i+1}\ot\dots \ot h_n\ot m,~~~
0\le i< n,~~~ ~~~
\nonumber \\ &&
\delta_n(h_0\ot\dots\ot h_n\ot m)= h_n^{(0)}h_0\ot h_1\dots
\ot h_{n-1}\ot h_{n}^{(1)}m,
\nonumber \\  &&
\sigma_i( h_0\ot\dots\ot h_n \ot m)= h_0 \ot \dots \ot h_i\ot
 1 \ot \dots \ot h_n\ot m,~~~
0\le i\le n,~~~ ~~~
\nonumber \\ &&\
\tau_n(h_0\ot\dots\ot h_n\ot m)= h_n^{(0)}\ot
  h_0\ot\dots \ot h_{n-1}\ot h_n^{(1)}m.  \nonumber
\eea
The homology of the  above cyclic module is by definition the Hopf cyclic homology of the Hopf algebra $\mathcal{H}$ with coefficients in the SAYD module $M$ and will be denoted by $\widetilde{HC}_{*}(\mathcal{H},M)$.
\ble\label{abcd} Let $\mathcal{H}$ and $\mathcal{K}$ be two Hopf algebras and $M$ and $N$  SAYD modules over $\mathcal{H}$ and $\mathcal{K}$ respectively. The following map is an isomorphism of cyclic modules
 \[\Omega_{n}: (\mathcal{H}\otimes \mathcal{K})^{\otimes n+1}\Box_{\mathcal{H}\otimes \mathcal{K}}(M\otimes N)\longrightarrow (\mathcal{H}^{\otimes  n+1}\Box_{\mathcal{H}}M)\otimes (\mathcal{K}^{\otimes n+1}\Box_{\mathcal{K}}N),\label{shuffle}\] given  by $$ ((h_{0}\otimes k_{0})\otimes\dots\otimes(h_{n}\otimes k_{n})\otimes(m\otimes r))\longmapsto (h_{0}\otimes \dots\otimes h_{n}\otimes m)\otimes(k_{0}\otimes \dots\otimes k_{n}\otimes r).$$

\ele
\bpf We prove $\Omega_{n}$ commutes with $\delta_{i}$ and $\tau_{n},$  where $ 0\leq i < n$. One can easily verify this for $\delta_{n}$ and degeneracies.
\begin{eqnarray}
&&(\delta_{i}\otimes \delta_{i})\Omega_{n}((h_{0}\otimes k_{0})\otimes\dots\otimes (h_{n}\otimes k_{n})\otimes (m\otimes r))\nonumber \\
&=& (\delta_{i}\otimes \delta_{i})((h_{0}\otimes\dots\otimes h_{n}\otimes m),(k_{0}\otimes \dots\otimes k_{n}\otimes r))\nonumber \\
&=&\delta_{i}(h_{0}\otimes\dots\otimes h_{n}\otimes m)\otimes \delta_{i}(k_{0}\otimes \dots\otimes k_{n}\otimes r)\nonumber \\
&=&(h_{0}\otimes\dots h_{i}h_{i+1}\otimes\dots\otimes h_{n}\otimes m)\otimes (k_{0}\otimes \dots k_{i}k_{i+1}\otimes\dots\otimes k_{n}\otimes r)\nonumber \\
&=& \Omega_{n}((h_{0}\otimes k_{0})\otimes \dots\otimes (h_{i}h_{i+1}\otimes k_{i}k_{i+1})\otimes \dots \otimes((h_{n}\otimes k_{n})\nonumber \\
&=& \Omega_{n}\delta_{i}((h_{0}\otimes k_{0})\otimes \dots\otimes (h_{n}\otimes k_{n})\otimes (m\otimes r)).\nonumber \end{eqnarray}

Since $\mathcal{H}$ is an $\mathcal{H}$-comodule algebra by comultiplication, we have:
\begin{eqnarray}
&&(\tau_{n}\otimes \tau_{n})\Omega_{n}((h_{0}\otimes k_{0})\otimes\dots\otimes(h_{n}\otimes k_{n})\otimes (m\otimes r))\nonumber \\
&=&\tau_{n}(h_{0}\otimes \dots\otimes h_{n}\otimes m)\otimes \tau_{n}(k_{0}\otimes \dots\otimes k_{n}\otimes r)\nonumber \\
&=& (h_{n}^{(1)}\otimes h_{0}\otimes \dots\otimes h_{n-1}\otimes h_{n}^{(2)}m)\otimes (k_{n}^{(1)}\otimes k_{0}\otimes \dots\otimes k_{n-1}\otimes k_{n}^{(2)}r).\nonumber \\
&=& \Omega_{n}(h_{n}^{(1)}\otimes k_{n}^{(1)})\otimes (h_{0}\otimes k_{0})\otimes \dots\otimes (h_{n-1}\otimes h_{n-1})\otimes (h_{n}^{(2)}m\otimes k_{n}^{(2)}r)\nonumber \\
&=& \Omega_{n}(h_{n}^{(1)}\otimes k_{n}^{(1)})\otimes (h_{0}\otimes k_{0})\otimes \dots\otimes (h_{n-1}\otimes h_{n-1})\otimes (h_{n}^{(2)}m\otimes k_{n}^{(2)}r)\nonumber \\
&=& \Omega_n\tau_{n}((h_{0}\otimes k_{0})\otimes \dots (h_{n}\otimes k_{n})\otimes (m\otimes r)).\nonumber
\end{eqnarray}
\epf

Using the previous lemma, the Eilenberg-Zilber isomorphism,  and the K\"{u}nneth formula we obtain the following proposition:
\bpr Let $\mathcal{H}$ and $\mathcal{K}$ be two Hopf algebras and $M$ and $N$ be SAYD modules over $\mathcal{H}$ and $\mathcal{K},$ respectively. We have the following isomorphism for Hopf Hochschild homology with coefficients

$$\widetilde{HH}_{n}(\mathcal{H}\otimes \mathcal{K},M\otimes N)\simeq \bigoplus_{{i+j=n}}\widetilde{HH}_{i}(\mathcal{H},M)\otimes \widetilde{HH}_{j}(\mathcal{K},N).$$ Also one has the following long exact sequence for  Hopf cyclic homology with coefficients:
$$\begin{CD}
... \longrightarrow \widetilde{HC}_{n}(\mathcal{H}\otimes \mathcal{K},M\otimes N)@>\mathcal{I} >> \bigoplus_{{i+j=n}} \widetilde{HC}_{i}(\mathcal{H},M)\otimes \widetilde{HC}_{j}(\mathcal{K},N)
@>S\otimes id-id \otimes S >>  \end{CD}$$
\[\begin{CD}
\bigoplus_{{i+j=n-2}}\widetilde{HC}_{i}(\mathcal{H},M)\otimes \widetilde{HC}_{j}(\mathcal{K},N)@> \partial >>
 \widetilde{HC}_{n-1}(\mathcal{H}\otimes \mathcal{K},M\otimes N)\longrightarrow ....
 \end{CD}\label{apcup}\] Furthermore if $(\widetilde{HC}(\mathcal{H},M)[-2m],S)_{-2m}$ satisfies the Mittag-Leffler condition and $\widetilde{HP}_{*}(\mathcal{H},M)$ is a finite dimensional vector space, then we obtain the following isomorphism for  the periodic Hopf cyclic homology with coefficients:
 $$\widetilde{HP}_{0}(\mathcal{H}\otimes \mathcal{K}, M\otimes N)\simeq \widetilde{HP}_{0}(\mathcal{H}, M)\otimes \widetilde{HP}_{0}(\mathcal{K}, N)\oplus  \widetilde{HP}_{1}(\mathcal{H}, M)\otimes \widetilde{HP}_{1}(\mathcal{K}, N),$$ and similarly for the odd case.
\epr

\bpr Let $\mathcal{H}$ be a cocommutative Hopf algebra and $M$ a SAYD module over $\mathcal{H}$ equipped with a map $\psi:M\longrightarrow M\otimes M$  satisfying the condition:
\[\psi(h.m)= \Delta(h).\psi(m). \label{shart2}\] The following map $$\rho_{n}=\Omega_{n} \Phi_{n} : \widetilde{C}_{n}(\mathcal{H}, M)\longrightarrow (\widetilde{C}(\mathcal{H},M)\times \widetilde{C}(\mathcal{H},M))_{n},$$ is a map of cyclic modules where
$$ \Phi_{n}= \psi\otimes \Delta^{\otimes n+1}:\widetilde{C}_{n}(\mathcal{H}, M) \longrightarrow \widetilde{C}_{n}(\mathcal{H}\otimes \mathcal{H}, M\otimes M).$$
\epr
\bpf We only check the commutativity of $\rho_{n}$ with $\delta_{n}$. The condition \eqref{shart2} is equivalent to
$$ (h.m)_{(1)}\otimes (h.m)_{(2)}= h^{(1)}.m_{(1)}\otimes h^{(2)}.m_{(2)}. $$ By cocommutativity of $\mathcal{H}$ we have:
\begin{eqnarray}
&& (\delta_{n}\otimes \delta_{n})\Omega_{n}(\Delta^{\otimes n+1}\otimes\psi)(h_{0}\otimes...\otimes h_{n}\otimes m)\nonumber \\
&=& \delta_{n}(h_{0}^{(1)}\otimes \dots \otimes h_{n}^{(1)}\otimes m_{(1)})\otimes \delta_{n}(h_{0}^{(2)}\otimes \dots \otimes h_{n}^{(2)}\otimes m_{(2)})\nonumber \\
&=& (h_{n}^{(1)(1)}h_{0}^{(1)}\otimes h_{1}^{(1)}\otimes\dots h_{n-1}^{(1)}\otimes h_{n}^{(1)(2)}m_{(1)})\nonumber \\
&& \otimes  (h_{n}^{(2)(1)}h_{0}^{(2)}\otimes h_{1}^{(2)}\otimes\dots h_{n-1}^{(2)}\otimes h_{n}^{(2)(2)}m_{(2)})\nonumber \\
&=& \Omega_{n} ((h_{n}^{(1)(1)}h_{0}^{(1)}\otimes (h_{n}^{(2)(1)}h_{0}^{(2)})\otimes (h_{1}^{(1)}\otimes h_{1}^{(2)})\nonumber \\
&&\otimes\dots\otimes (h_{n-1}^{(1)}\otimes h_{n-1}^{(2)})\otimes (h_{n}^{(1)(2)}m_{(1)}\otimes h_{n}^{(2)(2)}m_{(2)})\nonumber \\
&=& \Omega_{n} ((h_{n}^{(1)(1)}h_{0}^{(1)}\otimes (h_{n}^{(1)(2)}h_{0}^{(2)})\otimes (h_{1}^{(1)}\otimes h_{1}^{(2)})\nonumber \\
&& \otimes\dots\otimes (h_{n-1}^{(1)}\otimes h_{n-1}^{(2)})\otimes (h_{n}^{(2)}m)_{(1)}\otimes (h_{n}^{(2)}m)_{(2)})\nonumber \\
&=& \Omega_{n} (\Delta^{\otimes n+1}\otimes\psi)(h_{n}^{(1)}h_{0}\otimes h_{1}\otimes\dots\otimes h_{n-1}\otimes h_{n}^{(2)}m)\nonumber \\
&=& \Omega_{n} ( \Delta^{\otimes n+1}\otimes \psi)\delta_{n}(h_{0}\otimes\dots\otimes h_{n}\otimes m).\nonumber
\end{eqnarray}

\epf

Now we are ready to define the desired coproducts.

\bpr Let $\mathcal{H}$ be a cocommutative Hopf algebra, $M$ a SAYD module over $\mathcal{H}$  equipped with a map $\psi:M\longrightarrow M\otimes M$ satisfying  \eqref{shart2}. The following maps define coproducts for $\widetilde{HH}_{\ast}(\mathcal{H},M)$, $\widetilde{HC}_{\ast}(\mathcal{H}, M)$ and $\widetilde{HP}_{\ast}(\mathcal{H}, M)$:

\[\sqcup = \mathfrak{I} (AW_{n} \Omega_{n} \Phi_{n})^{*} : \widetilde{HH}_{n}(\mathcal{H},M)\longrightarrow \bigoplus_{{p+q=n}} \widetilde{HH}_{p}(\mathcal{H},M)\otimes \widetilde{HH}_{q}(\mathcal{H},M),\label{copp}\] and

\[\sqcup= \mathcal{I} (\widetilde{AW}_{n}\overline{\Omega}_{n} \overline{\Phi}_{n})^{*}: \widetilde{HC}_{n}(\mathcal{H},M)\longrightarrow \bigoplus_{{p+q=n}} \widetilde{HC}_{p}(\mathcal{H},M)\otimes \widetilde{HC}_{q}(\mathcal{H},M),\label{app2}\]
 where $\overline{\Phi}_{n}= \oplus_{i\geq 0} \Phi_{n-2i}$ and $\overline{\Omega}_{n} = \oplus_{i\geq 0} \Omega_{n-2i}$. Also if $(\widetilde{HC}(H,M)[-2m],S)_{-2m}$ satisfies the Mittag-Leffler condition and $\widetilde{HP}_{*}(\mathcal{H},M)$ is a finite dimensional vector space, then

\[\sqcup =\nabla ( \widetilde{AW}_{n}\underline{\Omega}_{n}\underline{\Phi}_{n})^{*}: \widetilde{HP}_{1}(\mathcal{H},M)\longrightarrow \widetilde{HP}_{0}(\mathcal{H},M)\otimes \widetilde{HP}_{1}(\mathcal{H},M)\oplus\widetilde{HP}_{1}(\mathcal{H},M)\otimes \widetilde{HP}_{0}(\mathcal{H},M),\]
 where $\underline{\Phi}_{n}= \oplus_{i\geq 0} \Phi_{2i}$
and $\underline{\Omega}_{n} = \oplus_{i\geq 0} \Omega_{2i}$ and similarly for the even case.
\epr

\bpf These are the results of the fact that by Lemma \ref{abcd}, the morphisms $\Phi$, $\overline{\Phi}$ and $\underline{\Phi}$ are maps of $(b,B)$-mixed complexes.
\epf

Recall that the Hopf Hochschild and Hopf cyclic homology of a group algebra $\mathcal{H}=kG$ are computed in \cite{kr1} and are given by
\[\widetilde{HH}_{n}(kG)\cong H_{n}(G,k),\] and
\[\widetilde{HC}_{n}(kG)\overset{\theta}\cong\bigoplus_{{i\geq 0}} H_{n-2i}(G,k)\label{mmm},\] where on the right hand side group homologies of $G$ with trivial coefficients appear. The coproduct in group homology is induced by the map
\[\sqcup_{Gr}(g_{1},\dots,g_{n})= \sum_{{k=0}}^{n} (g_{1}\otimes \dots\otimes g_{k}) \otimes (g_{k+1}\otimes\dots \otimes g_{n}).\]
\ble The coproduct for Hopf Hochschild homology of a group algebra $kG$ coincides with the coproduct of group homology, i.e.,
$$\sqcup(g_{1}\otimes\dots\otimes g_{n})= \sqcup_{Gr}(g_{1},\dots,g_{n}).$$
\ele
\bth Under the isomorphism \eqref{mmm}, the coproduct \eqref{app2} for the Hopf cyclic homology of the group algebra $kG$ with trivial coefficients coincides with the coproduct of group homology . Equivalently, the following diagram commutes on the level of homology.

\[\begin{CD}
\widetilde{C_{n}}(kG) @>\theta_{n}>> C_{n}^{Gr}(G,k)\\
@VV\cup' V @VV\cup_{Gr} V\\
(\widetilde{C}(kG) \otimes \widetilde{C}(kG))_{n}@>\theta'>> (C^{Gr}(G,k)\otimes C^{Gr}(G,k))_{n}
\end{CD}\]\\

where $\cup'= \widetilde{AW}\Omega \Delta^{\otimes n}$, $\widetilde{C_{n}}(kG)= kG^{n}$ and $\theta'=\sum_{{i\geq 0}}\theta_{n-2i}\otimes \theta_{n-2i}$.

\ethe

\section{Cup products in  Hopf cyclic homology}

Dual to previous coproducts, in this section we  define cup products for Hopf  Hochschild and Hopf cyclic homology of a commutative Hopf algebra.

Let $\mathcal{H}$ and $\mathcal{K}$ be two Hopf algebras. The map
 \[\times =sh_{p,q}: \widetilde{C}_{p}(\mathcal{H})\otimes \widetilde{C}_{q}(\mathcal{K})\longrightarrow (\widetilde{C}(\mathcal{H})\times \widetilde{C}(\mathcal{K}))_{p+q},\label{cup}\]  commutes with Hochschild boundaries and therefore by composing it with the  map $\Omega$ given by \eqref{shuffle}, one obtains the following map

$$ ((h_{1}\otimes \dots\otimes h_{p})\otimes (k_{1}\otimes \dots\otimes k_{q}))\longmapsto \sum_{\sigma} sign(\sigma)((h_{\sigma^{-1}(1)}\otimes 1)\otimes\dots (h_{\sigma^{-1}(p)}\otimes 1)\otimes (1 \otimes k_{\sigma^{-1}(1)}\otimes\dots (1\otimes k_{\sigma^{-1}(q))}),$$ where $\sigma$ runs through   $(p,q)$-shuffles. This  map induces a product  in  Hopf Hochschild homology:
$$\widetilde{HH}_{p}(\mathcal{H})\otimes \widetilde{HH}_{q}(\mathcal{K})\longrightarrow \widetilde{HH}_{p+q}(\mathcal{H}\otimes \mathcal{K}).$$

When $\mathcal{H}$ is commutative we can compose the above map with the the multiplication map $\mathcal{H}\otimes \mathcal{H}\longrightarrow \mathcal{H}$ to obtain a product on Hopf Hochschild homology.
\bpr Let $\mathcal{H}$ be a commutative Hopf algebra. The map
$$ \times: \widetilde{HH}_{p}(\mathcal{H})\otimes \widetilde{HH}_{q}(\mathcal{H})\longrightarrow \widetilde{HH}_{p+q}(\mathcal{H}),$$ induces  a structure of graded commutative algebra on $\widetilde{HH}_{\ast}(\mathcal{H})$.

\epr
\bpf Since $\mathcal{H}$ is commutative, the multiplication map $m : \mathcal{H}\otimes \mathcal{H}\longrightarrow \mathcal{H}$ is a Hopf algebra map. Now composing \eqref{cup}, for $\mathcal{H}=\mathcal{K}$, with the maps induced by $m$ and $\Omega$, provides us a product map:
$$\times: \widetilde{C}_{p}(\mathcal{H})\otimes \widetilde{C}_{q}(\mathcal{H})\longrightarrow \widetilde{C}_{p+q}(\mathcal{H}),$$
given by
$$(h_{1}\otimes \dots\otimes h_{p})\times (h_{p+1}\otimes \dots\otimes h_{p+q})=\sum_{\sigma} sign(\sigma)(h_{\sigma^{-1}(1)}\otimes\dots\otimes h_{\sigma^{-1}(p+q)}).$$ Here $\sigma$ runs over all $(p,q)$-shuffles. Therefore $\widetilde{C}_*(\mathcal{H})$ becomes a graded algebra.
\epf

 Now we define a cup product for Hopf cyclic homology of a commutative Hopf algebra.
 \bpr Let $\mathcal{H}$ and $\mathcal{K}$ be two Hopf algebras. The following map is a map of $(b,B)$-mixed complexes:
 $$\star: Tot\mathcal{B}\widetilde{C}_{p}(\mathcal{H})\otimes Tot\mathcal{B}\widetilde{C}_{q}(\mathcal{K})\longrightarrow Tot\mathcal{B}(\widetilde{C}(\mathcal{H})\times \widetilde{C}(\mathcal{K}))_{p+q+1},$$  given by
 \[ (x_{p}, x_{p-2}, \dots) \star (y_{p}, y_{p-2}, \dots)  = (Bx_{p}\times y_{q}, Bx_{p}\times y_{q-2}, \dots). \]
 \epr
 \bpf It is enough to show:
 \[b(Bx_{p} \times y_{q}) + B(Bx_{p}\times y_{q-2})= B(bx_{p}+ Bx_{p-2})\times y_{q} + (-1)^{p}Bx_{p}\times (by_{q}+ By_{q-2}),\label{rabt}\] and similarly for the other terms. Let $sh'_{i,j}(x,y)=x\times' y.$ One can verify  $Bx\times' y=0$ and $x\times'By=0$, for all $x$ and $y$ in the normalized complex. The equation  $[B, sh] + [b, sh']=0 $  is equivalent to
 $$B(x\times y)-(Bx\times y +(-1)^{i}x\times By)=-b(x\times'y)+bx\times'y+(-1)^{i}x\times' by.$$  By substituting $By$ instead of $y$ in the above equation we obtain:

   $$B(x\times By)= Bx \times By.$$ Now \eqref{rabt}  is the consequence of the fact that $[b, sh]= 0$.
 \epf

 \bth The map $\star$ induces the following associative product
 $$\widetilde{HC}_{p}(\mathcal{H})\otimes \widetilde{HC}_{q}(\mathcal{K}) \longrightarrow \widetilde{HC}_{p+q+1}(\mathcal{H}\otimes \mathcal{K}),$$ in Hopf cyclic homology. If we consider $\widetilde{HC}_{n}(\mathcal{H})$ of degree $n+1$, then the product is  graded commutative, i.e.,
 $$x \star y = (-1)^{(p+1)(q+1)}(y\star x) , $$ for $x \in \widetilde{HC}_{p}(\mathcal{H})$ and $y \in \widetilde{HC}_{q}(\mathcal{K}).$
 \ethe

 Now we are ready to define a cup product for Hopf cyclic homology of  commutative Hopf algebras:
 \bpr Let $\mathcal{H}$ be a commutative Hopf algebra. The  product $\star$ induces a graded commutative algebra structure on Hopf cyclic homology:
 $$\widetilde{HC}_{p}(\mathcal{H})\otimes \widetilde{HC}_{q}(\mathcal{H})\longrightarrow \widetilde{HC}_{p+q+1}(\mathcal{H}).$$
 \epr
 \bpf This is the consequence of that fact that when $\mathcal{H}$ is commutative, the product map $m$ is a Hopf algebra map.
 \epf

 One can see that
 \bpr The boundary map $\partial$ in the K\"{u}nneth long exact sequence \eqref{apcup} is the same as the product $\star$ in Hopf cyclic homology,
 $$\partial (x\otimes y)=x\star y,$$ where $x \in \widetilde{HC}_{p}(\mathcal{H})$ and $y \in \widetilde{HC}_{q}(\mathcal{K})$.
 \epr


\end{document}